\documentclass[final]{siamltex}

\usepackage[english]{babel}
\usepackage{amsmath}
\usepackage{amssymb}
\usepackage{epsfig}
\usepackage{array}
\usepackage{ifthen}
\usepackage{algorithm}
\usepackage{algorithmic}
\usepackage{subfigure}
\usepackage{booktabs}
\usepackage{aicescover}

\newcommand{\R}{\mathbb{R}}
\newcommand{\Rn}{\R^{n}}
\newcommand{\Rnn}{\R^{n \times n}}

\newcommand\norm[1]{\lVert#1\rVert}
\newcommand\order[1]{\mathcal{O}(#1)}
\newcommand\abs[1]{\lvert#1\rvert}

\newcommand{\binaryx}{{\it binary\_\hspace*{1pt}x}}
\newcommand{\binaryy}{{\it binary\_\hspace*{0pt}y}}

\begin{document}

\title{Improved Accuracy and Parallelism for MRRR-based Eigensolvers -- A Mixed Precision Approach\thanks{Financial support  
from the Deutsche Forschungsgemeinschaft (German Research Association)
through grant GSC 111 is gratefully acknowledged. Enrique S. Quintana-Ort\'{\i} was supported by project TIN2011-23283 and FEDER.}} 

\author{
  M.~Petschow\thanks{Aachen Institute for advanced study in Computational 
Engineering Science, RWTH Aachen University, Germany. Electronic address:
{\tt \{petschow,pauldj\}@aices.rwth-aachen.de}}
  \and E.~S. Quintana-Ort\'{i}\thanks{Depto. de Ingenier\'{\i}a y Ciencia de
    Computadores, Universidad Jaume I, 12071 Castell\'on, Spain. Electronic
    address: {\tt quintana@icc.uji.es}}
  \and P.~Bientinesi\footnotemark[2]
}

% Financial support from the Deutsche Forschungsgemeinschaft (German Research
% Association) through grant GSC 111 is gratefully acknowledged.
% Enrique S. Quintana-Ort\'{\i} was supported by project TIN2011-23283 and FEDER.

% \author{M. Petschow\,\footnote{RWTH Aachen, Aachen Institute for Advanced Study in
%   Computational Engineering Science, 52062 Aachen, Germany. Electronic
%   address: {\tt \{petschow,pauldj\}@aices.rwth-aachen.de}}\;,
% E.~S. Quintana-Ort\'{i}\,\footnote{Depto. de Ingenier\'{\i}a y Ciencia de Computadores,
%                             Universidad Jaume I,
%                             12071 Castell\'on, Spain. Electronic address:
%                             {\tt quintana@icc.uji.es}}\;, P. Bientinesi\,$^*$}

\aicescoverpage

\maketitle

\begin{abstract}
The real symmetric tridiagonal eigenproblem is of outstanding 
importance in numerical computations; it arises frequently as 
part of eigensolvers for standard and generalized 
dense Hermitian eigenproblems that are based on a
reduction to tridiagonal form.
For its solution, the algorithm of {\it Multiple Relatively
  Robust Representations} (MRRR) is among the fastest methods. 
Although fast, the solvers based on MRRR do not
deliver the same accuracy as competing methods like Divide \& Conquer or the
QR algorithm.   
In this paper, we demonstrate that the use of mixed precisions leads to
improved accuracy of MRRR-based eigensolvers {\it with limited or no
  performance penalty}. As a result, we obtain eigensolvers that are not
only equally or more accurate than the best available methods, but also --
under most circumstances -- faster and more scalable than the competition.
\end{abstract}

% to be changed
\begin{keywords} 
  symmetric eigenproblems, tridiagonal matrices, high-performance, MRRR
\end{keywords}

% to be changed
\begin{AMS}
65F15, 65Y05, 68W10
\end{AMS}

\pagestyle{myheadings}
\thispagestyle{plain}
\markboth{M.~Petschow, E.~S. Quintana-Ort\'{\i}, and
  P.~Bientinesi}{Improved accuracy for MRRR-based eigensolvers}

\section{Introduction}
\label{sec:intro}

In~\cite{Li02design}, the authors
describe how in libraries the use of ``higher internal
precision and mixed input/output types and precisions permits
[...] to implement some algorithms that are simpler, more accurate, and
sometimes faster.'' In particular, the internal use of higher precision
provides the 
library developer with extra precision and a wider range of values, which may
benefit the accuracy and robustness of numerical routines. 
In sharp contrast to software that uses
arbitrary precision to obtain any desired accuracy, the use of higher precision should not
lower performance  significantly if at all. In this paper, we employ
mixed precisions to improve not only the accuracy, but also the robustness and
scalability of eigensolvers based on the algorithm of {\it Multiple Relatively
  Robust Representations} (MRRR or MR$^3$ for
short)~\cite{DesignMRRR,Bientinesi:2005:PMR3,Vomel:2010:ScaLAPACKsMRRR,mr3smp,EleMRRR,Willems:Diss}.  

Direct methods for standard and generalized Hermitian eigenproblems often rely on a
reduction to real symmetric tridiagonal form~\cite{Bai:2000:TSA:357352}. Once the
problem is transformed, the {\it real symmetric tridiagonal eigenproblem} (STEP)
is the following: Given a tridiagonal matrix $T \in \Rnn$ (with $T =
T^*$, where $T^*$ denotes the transpose of $T$), find quantities $\lambda
\in \R$ and nonzero $z \in \Rn$ such that the equation
\begin{equation*}
T z = \lambda z 
\end{equation*}
holds. Without loss of generality, we assume $\norm{z} = 1$ hereafter,
where $\norm{\bullet}$ denotes the 2-norm. % Euclidean norm. 
For such a 
solution, $\lambda$ is called an \textit{eigenvalue} (of $T$) and $z$ an associated
\textit{eigenvector}; an eigenvalue together with an associated
eigenvector are said to form an \textit{eigenpair}, $(\lambda, z)$.
The Spectral Theorem~\cite{Parlett:1998:SEP} ensures the existence of
$n$ eigenpairs $(\lambda_i,z_i)$, $i \in
\{1,2,\ldots,n\}$, such that the
eigenvectors form a complete orthonormal set; that is, for all $i,j \in
\{1,2,\ldots,n\}$,   
\begin{equation*}
  z_j^* z_i = \left\{ 
\begin{array}{l l}
    1 & \quad \text{if $j = i$} \,, \\
    0 & \quad \text{if $j \neq i$} \,. \\
  \end{array} \right.
\end{equation*}
Since all eigenvalues are
real, they can be ordered:
\begin{equation*}
  \lambda_1 \leq \lambda_2 \leq \ldots \leq \lambda_i\leq \ldots \leq \lambda_n  \,,
\end{equation*}
where $\lambda_i$ is the $i$-th smallest eigenvalue of $T$. In
this paper, whenever the underlying matrix is not clear, we write
$\lambda_i[T]$ explicitly. The set of all eigenvalues of $T$ is denoted
$spec[T]$ and the spectral diameter is defined as $spdiam[T] = \lambda_n -
\lambda_1$. 
For a given 
index set $\mathcal{I}
\subseteq \{1,2,\ldots,n\}$,  
\begin{equation*}
  \mathcal{Z}_{\mathcal{I}} = \mbox{span} \{ z_i : i \in \mathcal{I}\}
\end{equation*}
denotes the invariant subspace associated with $\mathcal{I}$. As with the
eigenvalues, whenever the underlying matrix is not understood from context, we
write $\mathcal{Z}_{\mathcal{I}}[T]$ explicitly.

In many applications, only a subset of eigenpairs are of interest and need
to be computed. 
For the computed eigenpairs,
$(\hat{\lambda}_i, \hat{z}_i)$, $\norm{\hat{z}_i} = 1$ and $i \in \mathcal{I}
\subseteq \{1,2,\ldots,n\}$, the accuracy of the
results can be quantified by 
the {\it largest residual norm} and the {\it orthogonality}, 
respectively defined as 
\begin{equation}
  R = \max_{i \in \mathcal{I}} \frac{\| T \hat{z}_i - \hat{\lambda}_i \hat{z}_i
    \|_1}{\| T \|_1} \quad \mbox{and} \quad
  O = \max_{i \in \mathcal{I}} \max_{\substack{ j \in \mathcal{I} \\ j \neq i }} | \hat{z}_j^* \hat{z}_i | \label{def:defresortho}
  \,.
\end{equation}

A number of excellent algorithms for the STEP have been discovered. Among them,
Bisection and Inverse Iteration (BI)~\cite{Dhillon98currentinverse,Ipsen:1997:Invit},
the QR algorithm (QR)~\cite{qr61a,qr61b},  
Divide \& Conquer (DC)~\cite{dc81,dc94,dc95}, and the focus of this study,
MRRR~\cite{Dhillon:Diss,Dhillon:2004:Ortvecs,Dhillon:2004:MRRR,Fernando97,Parlett2000121,Willems:Diss}.  
These methods differ in various aspects: the number of floating point operations (flops)
they perform, the flop-rate at 
which the operations are executed, the amount of memory required, the
possibility of computing subsets 
of eigenpairs at reduced cost, the attainable accuracy, the simplicity and
robustness of the code, and the suitability for parallel computations.
Thus, the ``best'' algorithm is influenced by factors such as the problem
(e.g., dimension, subset, spectral distribution), the 
architecture (e.g., cache sizes, parallelism), external
libraries (e.g., Basic Linear Algebra Subprograms), and the specific
implementation of the algorithm (e.g., thresholds, optimizations). 

Demmel et al.~\cite{perf09} provide a detailed study of the performance and
accuracy of LAPACK's~\cite{laug} 
implementations of these four methods on various
architectures. They conclude that
($i$)~DC and QR are the most accurate
algorithms; ($ii$)~DC requires $\order{n^2}$ additional memory and therefore
much more than all 
the other algorithms;\footnote{Here and in the following, we use the notation
  $\order{x}$ informally as ``of the order of $x$ in magnitude.'' The notion
is used to hide moderate constants that are of no particular interest to our
discussion.} ($iii$)~DC and MRRR are {\em much faster} than QR and BI;
despite the fact that MRRR uses the fewest flops, DC is faster on certain classes of
matrices. If the full eigendecomposition is desired, DC is generally the method of
choice, but whether DC or MRRR is faster depends on the spectral distribution
of the input matrix; and ($iv$)~if only a subset of eigenpairs is desired, 
MRRR is the method of choice. 
The study is limited to {\it sequential} executions and does not take
into account the degree of parallelism the algorithms provide. However, various
studies~\cite{Bientinesi:2005:PMR3,Vomel:2010:ScaLAPACKsMRRR,EleMRRR,Tisseur:1999:PDC,mr3smp} of the performance and 
accuracy of parallel implementations come to similar conclusions.

To summarize, if all eigenpairs are computed, depending on the spectral
distribution of the input matrix, either DC or MRRR is the fastest method. If only
a subset of eigenpairs is desired, MRRR is the method of
choice. Unfortunately, MRRR delivers generally the least accurate results.   
These observations carry over to direct methods for the dense
eigenproblem based on a reduction to tridiagonal form. 
It is natural to ask whether the accuracy of the MRRR-based routines can be
improved to levels of other methods like QR or DC. 
Unfortunately, a general analysis of any MRRR implementation shows that, even if all 
requirements of the algorithm are fulfilled, one needs to expect
orthogonality of $\order{1000 n \varepsilon}$ -- with unit roundoff 
$\varepsilon$~\cite{Willems:framework}. Methods like QR and DC, however,
attain superior results with orthogonality of $\order{\varepsilon
  \sqrt{n}}$. In this paper, we present a practical solution that improves the
accuracy of MRRR. As a result, it becomes equally or more accurate than QR and DC.

Our solution resorts to the use of higher precision arithmetic. The motivation is
twofold: ($i$) MRRR is frequently the fastest algorithm and it might be
possible to trade (some of) its performance to obtain higher accuracy; ($ii$)
often MRRR is used in the context of direct methods for Hermitian
eigenproblems. While the tridiagonal stage is responsible for much of 
the ``loss'' of orthogonality in the final result, it has a lower complexity
than the reduction to
tridiagonal form. Thus, {\it even if it is necessary to spend more time in the
tridiagonal stage to improve accuracy, for sufficiently large matrices, the
overall run time will not be affected 
significantly}. As MRRR does not make use of
any level-3 Basic Linear Algebra Subprograms (BLAS)~\cite{blas90}, we do not
require in our mixed precision approach any optimized BLAS library for 
high precision, which might not be available. 

For any MRRR solver, we present how the use of mixed
precisions leads to more accurate results at very little or even no extra
costs in terms of performance. 
As a consequence, {\it MRRR is not only one 
of the fastest methods, but also becomes as accurate or even more accurate than the
competition}. Moreover, for direct methods based on a reduction to
tridiagonal form and MRRR, the tridiagonal eigensolver is responsible for
the inferior orthogonality compared with other methods. These solvers
benefit directly 
from our approach. 

\subsection{Related work and outline}
\label{sec:relatedwork}

The term {\it mixed precision algorithm} is sometimes synonymously
used for the following procedure: First solve the problem using a fast
low-precision arithmetic, and then refine the result to high accuracy using
a high-precision arithmetic.
This {\it mixed precision iterative refinement} approach
exploits the fact 
that  
there might exist a low-precision arithmetic faster
than that of the input/output data format. The larger the performance gap
between the two arithmetic, the more beneficial the
approach. Iterative refinement (with and without using mixed precisions) has
been most extensively studied for the 
solution of linear systems of equations~\cite{HIGHAM01101997},
but other 
operations such as the solution of Lyapunov equations   
also benefit from it~\cite{Benner:2011:MAS:2010586.2010622}.
 
We use of the term {\it mixed precision} in its more general form; that is,
using two or more different precisions for solving a problem. In particular,
we use a
higher precision in the more sensitive parts of an 
algorithm to obtain accuracy, which otherwise could not be achieved. This
approach is especially effective if the sensitive 
portion of the algorithm and/or the performance gap
between the two arithmetic is small. Similarly to the mixed
precision ideas for iterative refinement, the approach is quite
general and we believe it can benefit computations in numerous areas. 

The rest of the paper is organized as follows: In
Section~\ref{sec:mrrralgorithm}, we present the MRRR algorithm and its
accuracy limitations. The section mainly serves as a vehicle to introduce the
factors that influence accuracy, which are summarized in
Theorem~\ref{resthm}. The derivation of
the error bounds, which can be found in~\cite{Willems:framework}, is not important
for the understanding of our discussion. In
Section~\ref{sec:mixedgeneral}, we detail our mixed 
precision approach in a general setting. Besides presenting a way to improve
accuracy, we investigate the effects on memory
usage, robustness, and scalability.
In Section~\ref{sec:implementation}, we comment on an actual implementation
of our approach and 
elaborate on a number of practical issues. Finally, we present experimental
results of our mixed precision solvers in Section~\ref{sec:experiments}.

\section{The MRRR algorithm}
\label{sec:mrrralgorithm}

In this section, we present the MRRR algorithm to the detail necessary for the
later discussion. 
Our exposition is largely based
on~\cite{Willems:framework,Willems:twisted,Dhillon:2004:MRRR} and a
high-level description of the method given in Algorithm~\ref{alg:mrrr}. We
comment more thoroughly on various parts of the computation in the following.  
The goal is to present
Theorem~\ref{resthm} -- i.e., the factors that influence the accuracy of any
implementation of MRRR.

\paragraph{Preprocessing} Algorithm~\ref{alg:mrrr} assumes that the necessary
preprocessing is already performed; this includes the {\it scaling} of the entries, and
the so called {\it splitting} of the input matrix into principal submatrices
if off-diagonal entries are sufficiently small in
magnitude~\cite{Parlett:1998:SEP}. 
In this section, without any loss is generality, we assume that
the input matrix is (numerically) {\it irreducible}, i.e., 
no off-diagonal entry is ``small enough'' in magnitude that warrants setting it
to zero. The exact criterion is specified later.
\begin{algorithm}[ht]
 \small
    {\bf Input:} Irreducible symmetric tridiagonal $T \in \mathbb{R}^{n
    \times n}$; index set $\mathcal{I}_{in}  \subseteq \{1,\ldots,n\}$. \\
    {\bf Output:} Eigenpairs $(\hat{\lambda}_i, \hat{z}_i)$ with
    $i \in \mathcal{I}_{in}$.
    
    \vspace{1mm}

  \algsetup{indent=2em}
  \begin{algorithmic}[1]
    \STATE Select shift $\mu \in
    \mathbb{R}$ and compute $M_{root} = T - \mu I$. \label{line:mrrr:root} 
    \STATE Perturb $M_{root}$ by a ``random'' relative amount bounded by a
    small multiple of $\varepsilon$. \label{line:mrrr:perturb}
    \STATE Compute $\hat{\lambda}_i[M_{root}]$ with $i \in \mathcal{I}_{in}$
    to relative accuracy sufficient for classification. 
    \label{line:mrrr:initialeigvals} 
    \STATE Form a work queue $Q$ and enqueue task
    $\{M_{root}, \mathcal{I}_{in}, \mu\}$. \label{line:mrrr:enddlarre}
    \WHILE{$Q$ not empty} 
    \STATE Dequeue a task $\{M, \mathcal{I}, \sigma\}$. \label{line:mrrr:dequeue}
    \STATE Partition $\mathcal{I} = \bigcup_{r=1}^R \mathcal{I}_r$  according 
    to the separation of the eigenvalues. \label{line:mrrr:initialpartitioning}
    \FOR{$r=1$ {\bf to} $R$}
    \IF{$\mathcal{I}_r = \{i\}$} \label{line:mrrr:ifstatement} 
    \STATE // {\it process well-separated eigenvalue associated with
      singleton $\mathcal{I}_r$} //
    \STATE Perform Rayleigh quotient iteration (guarded by bisection) to
    obtain eigenpair $(\hat{\lambda}_i[M],\hat{z}_i)$ with sufficiently
    small residual norm, $\norm{M \hat{z}_i - \hat{\lambda}_i[M]
      \,\hat{z}_i}$. \label{line:mrrr:rqi}
    \STATE Return $\hat{\lambda}_i[T] = \hat{\lambda}_i[M] + \sigma$ and
     $\hat{z}_i$. \label{line:mrrr:returneigenpair} 
    \ELSE
    \STATE // {\it process cluster associated with $\mathcal{I}_r$} //
    \STATE Select shift $\tau \in
    \mathbb{R}$ and compute $M_{shifted} = M - \tau
    I$. \label{line:mrrr:shifting}  
    \STATE Refine $\hat{\lambda}_i[M_{shifted}]$ with $i \in \mathcal{I}_r$ to
    sufficient relative 
    accuracy. \label{line:mrrr:refine} 
    \STATE Enqueue $\{M_{shifted}, \mathcal{I}_r, \sigma +
    \tau\}$. \label{line:mrrr:partition} 
    \ENDIF
    \ENDFOR
    \ENDWHILE \label{line:mrrr:end}
  \end{algorithmic}
  \caption{\ MRRR}
  \label{alg:mrrr}
\end{algorithm}

\paragraph{Choice of representations}
In order for Algorithm~\ref{alg:mrrr} 
to work, the representation of tridiagonals (i.e., $M_{root}$ and
$M_{shifted}$) by their diagonal and off-diagonal entries must be abandoned
and alternative representations must be used. Any $2n-1$ (or less) scalars
together with a mapping that define the entries of a symmetric tridiagonal
is called a {\it representation}~\cite{Willems:framework}. We distinguish 
between the data of the representation, which are floating point numbers,
and the underlying tridiagonal, which is generally not exactly representable
 in the same finite precision format. There are multiple candidates
-- existence assumed -- for providing representations of tridiagonals:

\begin{enumerate}
\item {\it Lower bidiagonal factorizations} of the form $T = L D L^*$, and {\it
    upper bidiagonal factorizations} of the form $T = U \Omega U^*$, where $D = 
  \mbox{diag}(d_1,d_2,\ldots,d_n) \in \Rnn$ and $\Omega = 
  \mbox{diag}(\omega_1,\omega_2,\ldots,\omega_n) \in \Rnn$ are diagonal,  
  $L \in \Rnn$ and $U \in \Rnn$ are respectively unit lower bidiagonal and
  unit upper bidiagonal. 
\item A generalization of the above are the so called {\it twisted
    factorizations} or {\it
    BABE-factorizations}~\cite{Fernando97}, $T = N_k 
  \Delta_k N_k^*$, where 
  $k$ denotes the {\it twist index}. The $k \times k$ leading principle
  submatrix of $N_k \in \Rnn$ is unit lower bidiagonal (determined by the non-trivial
  entries $\ell_1,\ldots,\ell_{k-1}$), and the $(n-k+1) \times (n-k+1)$
  trailing principle submatrix of $N_k$ is unit upper bidiagonal (determined by the non-trivial
  entries $u_k,\ldots,u_{n-1}$); the matrix $\Delta_k =
\mbox{diag}(d_1,\ldots,d_{k-1},\gamma_k,\omega_{k+1},\ldots,\omega_n) \in \Rnn$ 
is diagonal. 
% Twisted factorizations for all indices $1 \leq k \leq n$ are almost entirely
% defined by elements of the bidiagonal factorizations $T = L D L^* = U
% \Omega U^*$, and only $\gamma_k$ has to be computed. 
Although 
it was known that these factorizations can additionally serve as
representations of the intermediate
matrices~\cite{Dhillon:Diss,Dhillon:2004:Ortvecs},
their benefits were only demonstrated 
recently~\cite{Willems:Diss,Willems:twisted}. 
Due to their additional degree of freedom in choosing $k$, the twisted
factorizations
% , which includes the bidiagonal factorizations as special
% cases, 
are superior to lower or upper bidiagonal factorizations. 
Besides representing intermediate tridiagonals, twisted factorizations are
essential in computing accurate eigenvectors~\cite{Fernando97,Dhillon:2004:Ortvecs}.
\item {\it Blocked
    factorizations} are further generalizations of bidiagonal and twisted factorizations. The
  quantities $D$, $\Omega$, and  
  $\Delta_k$ are block diagonal, with blocks of size $1 \times
  1$ or $2 \times 2$. The other factors -- $L$, $U$, and $N_k$ -- are
  partitioned conformally, with one or the $2 \times 2$ identity as diagonals. These types of
  factorizations contain the unblocked bidiagonal and twisted factorizations
  as special cases. With their great flexibility, the blocked factorizations
  have been used very successfully within the MRRR
  algorithm~\cite{Willems:Diss,Willems:blocked}. 
\end{enumerate}
All these factorizations are determined by $2n-1$ scalars, the data. 
For instance, for lower bidiagonal factorizations, the $2n-1$ floating
point numbers $d_1,\ldots,d_n,\ell_1,\ldots,\ell_{n-1}$ determine a
tridiagonal; such representation by the non-trivial entries of the
factorization is called an $N$-representation. 
Similarly, the floating point numbers
$d_1,\ldots,d_n,e_1,\ldots,e_{n-1}$ represent a tridiagonal -- with $e_i = d_i \ell_i$,
$1 \leq i \leq n-1$,
being $T$'s off-diagonal elements; such representation, including $T$'s
off-diagonal elements, is called an $e$-representation. Besides the $N$- and
the $e$-representation, the $Z$-representation is introduced
in~\cite{Willems:twisted} for bidiagonal and twisted factorizations. For
blocked factorizations, a variant of the $e$-representation is commonly
used~\cite{Willems:blocked}. 
Other quantities that are computed using the
(primary) data are called {\it secondary} or {\it derived} data. For
instance, $T$'s off-diagonal 
elements are secondary for an $N$-representation while
being primary for an $e$-representation. While the details are not relevant
for our discussion, it is important to note that there are different
variants to represent tridiagonals -- each one with slightly different properties.
Subsequently, we do not distinguish between the representation of a
tridiagonal and the tridiagonal itself; that is, it is always implied that
tridiagonals are represented in one of the above forms.

\paragraph{The representation tree}
The unfolding of Algorithm~\ref{alg:mrrr} is best described as a tree of
representations~\cite{Dhillon:Diss,Dhillon:2004:MRRR,Willems:framework}. Each task
$\{M, \mathcal{I}, \sigma\}$ or just $\{M, \mathcal{I}\}$ (Line~\ref{line:mrrr:dequeue} of
Algorithm~\ref{alg:mrrr}) is connected to a
node in the tree; that is, all nodes consist of a representation and an
index set. $\{M_{root},\mathcal{I}_{in}\}$ is associated with the root node (hence the name).
All other tasks are connected to ordinary nodes. Each node has a
depth associated with it: the number of edges on the unique path to it from
the root. The maximum depth for all nodes is denoted $d_{max}$. The edges
connecting internal nodes are associated with the spectrum shifts $\tau$ that are
performed in Line~\ref{line:mrrr:shifting} of Algorithm~\ref{alg:mrrr}. 

\paragraph{Factors influencing MRRR's accuracy} The analysis
in~\cite{Willems:framework} 
-- a streamlined version of the proofs
in~\cite{Dhillon:2004:MRRR,Dhillon:2004:Ortvecs} --  shows that, if suitable
representations are found, the 
computed eigenpairs enjoy a small residual norm and are mutually
(numerically) orthogonal. 

\begin{theorem}[Accuracy] 
Let $\hat{\lambda}_i[M_{root}]$ be computed (exactly)\footnote{The
  assumption can be removed; we simply stated the theorem as it can be found
in~\cite{Willems:Diss,Willems:framework}.} by applying the
spectrum shifts to eigenvalue $\hat{\lambda}_i[M]$ obtained by the Rayleigh quotient
iteration in Line~\ref{line:mrrr:rqi} of Algorithm~\ref{alg:mrrr}.  
Provided all the requirements, which are discussed below, are satisfied, it holds
\begin{equation*}
\norm{M_{root}\, \hat{z}_i - \hat{\lambda}_i[M_{root}]  \, \hat{z}_i} \leq
\left( \norm{r^{(local)}} + \gamma \, spdiam[M_{root}] \right)
\frac{1+\eta}{1-\eta}
\label{residualbound}
\end{equation*}
with $\gamma = k_{elg} n \left( d_{max}(\xi_{\,\downarrow} + \xi_{\uparrow})
  +  \alpha \right) + 2(d_{max} +1) \eta$. 
Furthermore, we have for any computed eigenvectors $\hat{z}_i$ and $\hat{z}_j$, $i \neq j$, 
\begin{equation*}
|\hat{z}_i^* \hat{z}_j| \leq 2 \left( \mathcal{R} n \varepsilon +  \frac{k_{rr}
  n (\xi_{\,\downarrow} + \xi_{\uparrow}) d_{max}}{gaptol} \right) \,.
\label{orthogonalitybound}
\end{equation*}
where $\mathcal{R} n \varepsilon = k_{rr} n \alpha / gaptol + k_{rs} n
\varepsilon/gaptol + \eta$.
\label{resthm}
\end{theorem}
A proof of the theorem can be found
in~\cite{Willems:Diss,Willems:framework}. 
Provided the representation of
$M_{root}$ is computed in a backward stable manner, a small residual norm with
respect to $M_{root}$ implies a small residual, 
$\order{n\varepsilon \norm{T}}$, with respect to the input matrix $T$. 

The rest of this section serves to convey the meaning of all the
parameters involved in the theorem. In Section~\ref{sec:mixedgeneral}, we
furthermore discuss their effects on performance, robustness,
and parallel scalability.  

\paragraph{Shifting the spectrum ($\xi_{\,\downarrow}, \xi_{\uparrow}$)}
The spectrum shifts of Line~\ref{line:mrrr:shifting} leave the eigenvectors 
unchanged in exact arithmetic; this invariance is lost in finite precision. An essential
ingredient of MRRR is the use of special forms of Rutishauser's {\it
  Quotienten-Differenzen} (qd) algorithm to perform the shifts. 
Given a representation for $M$, we require that the representation
for $M_{shifted} = M - \tau I$ is computed in an {\it element-wise mixed relative 
stable} way, i.e., $\widetilde{M}_{shifted} = \widetilde{M} - \tau I$ holds
exactly for small element-wise relative perturbations of the data for
$M_{shifted}$ and $M$. For all shifts performed in the algorithm, these
perturbations must be bounded by   
$\xi_{\uparrow} = \order{\varepsilon}$ and $\xi_{\,\downarrow} =
\order{\varepsilon}$, respectively. Algorithms that implement the spectrum
shifts for different forms of representations are presented
in~\cite{Dhillon:2004:Ortvecs,Willems:twisted,Willems:blocked}. 

\paragraph{Requirements on the representations ($k_{rr}, k_{elg}$)}
In order to ensure that the computed eigenpairs enjoy small residual norms
with respect to the input matrix and that the eigenvectors are
numerically orthogonal, the representations in
Line~\ref{line:mrrr:shifting}, $M_{shifted} = M - \tau I$, need to be chosen
with care. By selecting appropriate shifts $\tau$, representations that are
{\it relatively robust} and exhibit {\it conditional element growth} are
selected. Before we define the meaning of these two concepts, we give a
brief definition of a {\it relative gap}. 

\begin{definition}[Relative gap] Let $T \in \Rnn$ be an irreducible
  symmetric tridiagonal matrix with eigenvalues $\{\lambda_i: 1 \leq i
  \leq n\}$ and let
  $\mathcal{I} \subset \{ 1, 2, \ldots , n\}$ be an index set.  
  The relative gap connected to $\mathcal{I}$ is defined as
  \begin{equation*}
    \mbox{\em relgap} (\mathcal{I}) = \min \left\{ \frac{|\lambda_j -
    \lambda_i|}{|\lambda_i|} : i \in 
    \mathcal{I}, j \notin \mathcal{I} \right\}
  \end{equation*}
  where quantities $|\lambda_j -
  \lambda_i|/|\lambda_i|$ are $\infty$ if $\lambda_i = 0$.
\label{def:relgapindexset}
\end{definition}

\begin{definition}[Relative robustness] 
Let $T$ (given by any representation), $\lambda_i$ and $\mathcal{I}$ be as
in the previous definition. Furthermore, let $\mathcal{Z}_{\mathcal{I}}$ be
the invariant subspace associated with $\mathcal{I}$. We say 
  that the representation of $T$ is relatively robust for $\mathcal{I}$
  if for all element-wise relative perturbation in the data bounded by
  $\xi \ll 1$ and for all $i \in \mathcal{I}$, we have 
  \begin{align*}
    |\tilde{\lambda}_i - \lambda_i| &\leq k_{rr} n \xi |\lambda_i| \,,  \\
    \sin \angle (\tilde{\mathcal{Z}}_{\mathcal{I}},\mathcal{Z}_{\mathcal{I}}) &\leq \frac{k_{rr} n \xi}{\mbox{\em relgap}(\mathcal{I})} \,,
  \end{align*}
where $\tilde{\lambda}_i$ and $\tilde{\mathcal{Z}}_{\mathcal{I}}$ denote the
eigenvalues and the corresponding invariant subspaces of the perturbed matrices,
respectively; $\angle(\tilde{\mathcal{Z}}_{\mathcal{I}},\mathcal{Z}_{\mathcal{I}})$
denotes the largest principle angle; 
and $k_{rr}$ is moderate constant, say about 10.\footnote{According to
  \cite{Willems:Diss,Willems:framework}, the requirement on the eigenvalues
  can be removed entirely: for $\mathcal{I}=\{i\}$, by
  Theorem~\ref{thm:gapthm} stated below, the second condition implies 
  the first up to a small constant provided $gap(\tilde{\lambda}_i) \approx
  gap(\{i\})$.  Similarly, if $\mathcal{I} = \{p,\ldots,q\}$ is not a
  singleton, the second term implies that $\tilde{\lambda}_i \in [\lambda_p
  - k_{rr} n \xi |\lambda_p|, \lambda_q - k_{rr} n \xi |\lambda_q|]$ for all $i \in \mathcal{I}$. 
}
\label{def:RRR}
\end{definition}

\begin{definition}[Conditional element growth] 
A representation for a real symmetric tridiagonal, $M$, exhibits conditional
element growth with respect to the index set $\mathcal{I} \subset \{1,\ldots,n\}$,
if for any element-wise relative perturbation in the data
bounded by $\xi \ll 1$ (leading to a perturbed tridiagonal $\widetilde{M}$) and each $i \in 
\mathcal{I}$, it holds
\begin{align*}
  \norm{\widetilde{M} - M} &\leq \mbox{\em spdiam}[M_{root}] \,, \quad
  \mbox{and} \\ % spdiam(T) 
  \norm{(\widetilde{M} - M) \hat{z}_i} &\leq k_{elg} n \xi \cdot \mbox{\em spdiam}[M_{root}]
  \,, % spdiam(T) 
\end{align*}
where $\hat{z}_i$ are the computed eigenvectors, and $k_{elg}$ is a moderate
constant, say about 10. 
\end{definition}

In Line~\ref{line:mrrr:shifting} of Algorithm~\ref{alg:mrrr}, we need to
ensure that $M_{shifted}$ as well as $M$ are relatively robust for
$\mathcal{I}_{r}$ and that $M_{shifted}$ features conditional element growth
for $\mathcal{I}_r$. 
In this paper, we are not concerned how to ensure that the involved representations satisfy
the requirements; this is the topic
of~\cite{perturbLDL,Parlett2000121,Dhillon:2004:Ortvecs,Willems:framework}.\footnote{In
  particular, it is {\it not} necessary to compute the eigenvectors in order
  to give bounds on the conditional element growth.} 
We 
remark however that there exist the danger that no suitable representation
that passes the test for the requirements can be found. In this case,
it is common to select a promising representation, which might not fulfill the
requirements. As a consequence, the accuracy of Theorem~\ref{resthm} is not
guaranteed anymore.    

\paragraph{Classification of the eigenvalues ($gaptol$)} 

While the above requirements on the representations pose a restriction on the
choice of shifts $\mu$ and $\tau$ in
Lines~\ref{line:mrrr:root} and \ref{line:mrrr:shifting}, the main goal is to
chose shifts such that, in the next iteration, the partitioning $\mathcal{I} = \bigcup_{r=1}^R
\mathcal{I}_r$, splits the index set into at least two subsets so that progress
in the algorithm is guaranteed. 
The partitioning is done according to the 
separation of the eigenvalues and must ensure two requirements: For a given
tolerance $gaptol$, say $10^{-3}$, ($i$) $relgap(\mathcal{I}_r) \geq gaptol$ , and ($ii$)
whenever $\mathcal{I}_r = \{i\}$ is a singleton, $relgap(\hat{\lambda}_i) \geq 
gaptol$. The latter relative gap is thereby defined as 
\begin{equation*}
relgap(\hat{\lambda}_i) = \frac{gap(\hat{\lambda}_i)}{|\lambda_i|} \quad
\mbox{with} \quad gap(\hat{\lambda}_i) = \min_{j \neq i} \left\{
  |\hat{\lambda}_i - \lambda_j| \right\} \,.
\end{equation*} 

For all $i \in \mathcal{I}$, let $\hat{\lambda}_i$ denote the midpoint point
of a computed interval of uncertainty
$[\underline{\lambda}_i,\overline{\lambda}_i]$ containing eigenvalue
$\lambda_i$. To achieve the desired partitioning of $\mathcal{I}$, let $j,j+1
\in \mathcal{I}$ and define 
\begin{equation*}
  reldist(j,j+1) = \frac{\underline{\lambda}_{j+1} -
    \overline{\lambda}_j}{\max\{ |\underline{\lambda}_j| ,
    |\overline{\lambda}_j|, |\underline{\lambda}_{j+1}| ,
    |\overline{\lambda}_{j+1}| \}}
\end{equation*}
as a measure of the relative gap. If $reldist(j,j+1) \geq gaptol$, then $j$ and
$j+1$ belong to different subsets of the partition. Additionally, this criterion
based on the relative separation can be amended by a criterion based on the
absolute separation of the eigenvalues~\cite{VoemelRefinedTree2007tr}. After
partitioning, each index set 
$\mathcal{I}_r$ with $|\mathcal{I}_r| > 1$ is associated with a {\it cluster}
of eigenvalues, $\{\hat{\lambda}_i : i \in \mathcal{I}_r\}$. Similarly, each
singleton $\mathcal{I}_r = \{i\}$ is associated with a {\it
  well-separated} eigenvalue $\hat{\lambda}_i$.

In order to reliably classify the
eigenvalues, they should be approximated in
Lines~\ref{line:mrrr:initialeigvals} and \ref{line:mrrr:refine} to relative
accuracy of about $gaptol$: that is, at least 
\begin{equation}
|\hat{\lambda}_i - \lambda_i| \lesssim gaptol \cdot |\lambda_i| \,.
\label{eq:relaccgaptol}
\end{equation}
The above criterion can be relaxed for
eigenvalues with a large gap to the rest of the
spectrum~\cite{DesignMRRR}. Commonly, the 
eigenvalues are computed by some form of 
bisection~\cite{Parlett:1998:SEP}. In 
particular, in Line~\ref{line:mrrr:refine}, we already have good
approximations to the eigenvalues, which can be refined by bisection to the
desired accuracy. 

The parameter $gaptol$ is so important that it influences almost all parts
of the algorithm. Since the error bounds in Theorem~\ref{resthm} are
proportional to $1/gaptol$, the value indicates how much accuracy we are willing to lose in the
computation. For many applications, this limits the choice to values larger than about
$10^{-3}$. However, we cannot use values much larger than
$10^{-3}$ as otherwise it becomes impossible to make progress by breaking
clusters. 

As a side note: the condition $relgap(\mathcal{I}_r) \geq gaptol$, together with the mixed
relatively stable computation of the spectrum shifts, implies that the
associated invariant subspaces are not 
perturbed too much due to rounding
errors, i.e., $\sin \angle
(\mathcal{Z}_{\mathcal{I}_r}[M_{shifted}],\mathcal{Z}_{\mathcal{I}_r}[M]) \leq
k_{rr} n (\xi_{\,\downarrow} + \xi_{\uparrow}) / gaptol$. After shifting, we can
therefore hope to compute an orthonormal 
basis for such a subspace, which is {\it automatically} numerically
orthogonal to the subspace spanned by the other eigenvectors. This is the
main idea behind the MRRR algorithm.

\paragraph{Rayleigh quotient iteration ($k_{rs}, \alpha, \eta$)}

Finally, in Line~\ref{line:mrrr:rqi} of Algorithm~\ref{alg:mrrr}, eigenpairs of well-separated
eigenvalues are computed via the Rayleigh quotient iteration (RQI).  
Given an approximation $\hat{\lambda}_i[M]$ and a representation $M$ that
is relatively robust for $\{i\}$, a key ingredient of MRRR is the ability to
compute an accurate eigenvector approximation $\hat{z}_i$ such that $\sin
\angle(\hat{z}_i,z_i) = \order{n\varepsilon/gaptol}$; see~\cite{Dhillon:2004:Ortvecs} for a
proof. This is certainly achieved by driving the local residual 
norm below a specified threshold
\begin{equation}
\norm{r^{(local)}} = \norm{M \hat{z}_i -
  \hat{\lambda}_i[M]\,\hat{z}_i} \leq k_{rs}  \cdot
gap\left(\hat{\lambda}_i[M]\right) \cdot \frac{n \varepsilon}{gaptol} \,,
\label{localresbound}
\end{equation}
where $k_{rs}$ is $\order{1}$. 
In this case, the so called Gap Theorem gives the desired bound on the error
angle $\angle(\hat{z}_i,z_i)$. 
\begin{theorem}[Gap Theorem]
Given a symmetric matrix $T \in \Rnn$ and an approximation $(\hat{\lambda},
\hat{z})$, $\| \hat{z} \| = 1$, 
to the eigenpair $(\lambda,
z)$, with $\hat{\lambda}$ closer to $\lambda$ than to any other eigenvalue,
let $r$ 
be the residual $T \hat{z} -
\hat{\lambda} \hat{z}$; then
\begin{equation}
\sin \angle (\hat{z},z) \leq \frac{\| r \|}{gap(\hat{\lambda})} \,.
\label{eq:gapthmmain}
\end{equation}
The residual norm is minimized if $\hat{\lambda}$ is the Rayleigh quotient
of $\hat{z}$, $\hat{\lambda} = \hat{z}^* T \hat{z}$. In 
this case, 
\begin{equation}
\frac{\| r \|}{spdiam[T]} \leq \sin \angle (\hat{z},z)  \quad
\mbox{and} \quad |\hat{\lambda} - \lambda| \leq \min \left\{
 \norm{r}, \frac{\norm{r}^2}{gap(\hat{\lambda})} \right\} \,.
\label{eq:gapthmeq2}
 \end{equation}
\label{thm:gapthm}
\end{theorem}
A proof of the theorem can be found for instance in~\cite{Parlett:1998:SEP,Willems:Diss}. 

In general, a residual norm such as in~\eqref{localresbound} cannot be guaranteed; it is only
possible to show that it holds for a small element-wise relative
perturbation of the data of $M$ bounded by $\alpha$ and the computed
eigenvector $\hat{z}_i$ bounded by $\eta$ -- with $\alpha = \order{\varepsilon}$ and
$\eta = \order{n \varepsilon}$. For our purposes, this
detail is not important. 
Nonetheless, Theorem~\ref{resthm} takes
this fact into account. Note that even in rare cases where \eqref{localresbound} is not
fulfilled, the small error angle together with \eqref{eq:gapthmeq2} imply
$\norm{r^{(local)}} = \order{spdiam[M] \cdot n\varepsilon/gaptol}$.

In the RQI, the $j$-th iteration consists of four steps: ($i$) For all $1 \leq k \leq
n$, compute the twisted factorizations $N_k \Delta_k N_k = M -
\hat{\lambda}_i^{(j)} I$; ($ii$) determine $s = \operatorname*{arg\,min}_k
|\gamma_k|$, where $\gamma_k$ is the $k$-th element of $\Delta_k$ (see
above); ($iii$) solve the linear system $N_s \Delta_s N_s^* \,
\hat{z}_i^{(j)} = \gamma_s e_s$, which is equivalent to the system $N_s^* \,
\hat{z}_i^{(j)} = e_s$; ($iv$) use the 
Rayleigh quotient correction %(RQC) 
term to update the eigenvalue
$\hat{\lambda}_i^{(j+1)} = \hat{\lambda}_i^{(j)} +
\gamma_s/\norm{\hat{z}_i^{(j)}}^2$. The residual norm is
approximated by $|\gamma_s|/\norm{\hat{z}_i^{(j)}}$ and the process is stopped if
\eqref{localresbound} is satisfied. 
In order to always converge, the stopping criterion is amended and the
iteration stopped when $\hat{\lambda}_i$ is not improved anymore, i.e.,
$|\gamma_s|/\norm{\hat{z}_i^{(j)}}^2 = \order{\varepsilon
  |\hat{\lambda}_i|}$~\cite{NLA:NLA493}.  
An alternative approach to RQI is to refine the eigenvalue approximation to full
precision (i.e., $|\hat{\lambda}_i - \lambda_i| = \order{n \varepsilon
  |\lambda_i|}$), and then perform only a single step of RQI. This approach is used 
whenever RQI fails to converge to the correct eigenvalue~\cite{DesignMRRR}.

\section{Mixed precision MRRR}
\label{sec:mixedgeneral}

The exact values of the parameters in Theorem~\ref{resthm} differ
slightly for various implementations of the algorithm and need not to be
known exactly in the following analysis. The bounds on the residual norm and
orthogonality are {\it theoretical}. It is useful to translate what the bounds mean in practice: 
with reasonable parameters, realistic {\it practical} bounds on the residual
norm and on the orthogonality are $n \varepsilon$ and $1000 n \varepsilon$, respectively. 
In order to obtain accuracy similar to that of the best available methods,
we need to trade the 
dependence on $n$ by a dependence on $\sqrt{n}$. Furthermore, it is
necessary to reduce the orthogonality by about three orders of
magnitude.

\subsection{A solver using mixed precisions}

The technique is simple, yet powerful: Inside the algorithm, 
we use a precision higher than of the input/output in order to improve accuracy. 
A similar idea was already mentioned in~\cite{Dhillon:Diss},
in relation to a preliminary version of the MRRR algorithm,
but was never pursued further.
With many implementation and algorithmic advances since then (e.g.,
\cite{NLA:NLA493,Dhillon05gluedmatrices,Bientinesi:2005:PMR3,Willems:twisted,Willems:blocked}),  
it is appropriate to  
investigate the approach in detail.
To this end, we build a tridiagonal eigensolver that 
differentiates between two precisions: ($i$)~the  
input/output precision, say \binaryx, and ($ii$)~the working precision,
\binaryy, with $y \geq x$. If $y = x$, we have the original 
situation of a solver based on one precision; in this case, the following analysis
is easily adapted to situations in which we are satisfied with {\em less} accuracy
than achievable by MRRR in $x$-bit arithmetic. Since we are interested in accuracy that
cannot be accomplished in $x$-bit arithmetic, we
restrict ourselves to the case $y > x$. Provided the unit roundoff of the $y$-bit format is
sufficiently smaller than the unit roundoff of the $x$-bit format, say four
or five orders of magnitude, we show how to obtain, for practical matrix sizes, improved
accuracy to the desired level. 

Although any $x$-bit and $y$-bit floating point format might be chosen, in
practice, only those shown in
Table~\ref{tab:precisions} are used in high-performance libraries. 
For example, for a {\it binary32} input/output format (single precision), we
might use a {\it binary64} working format (double precision). Similarly, for
a {\it binary64} input/output format, we might use a {\it binary80} or {\it binary128} working
format (extended or quadruple precision). For these three configurations, we use the terms
{\it single/double}, {\it double/extended}, and {\it
  double/quadruple}. Practical issues for their implementation are  
discussed in Section~\ref{sec:implementation}. In this section, however, we
concentrate on the generic case of \binaryx/\binaryy. In general, when we
refer to \binaryx, we mean both the $x$-bit data type and its unit roundoff $\varepsilon_x$. 
\begin{table}[htb]
  \begin{center}
\begin{tabular}[thb]{l@{\quad}l@{\quad}l@{\quad}l@{\quad}l@{\quad}}
\toprule
{\bf Name} & {\bf IEEE-754} & {\bf Precision} & {\bf Support}  \\
\midrule
single    & binary32   & $\varepsilon_s = 2^{-24}$ &  Hardware       \\ 
double    & binary64   & $\varepsilon_d = 2^{-53}$ &  Hardware       \\ 
extended  & binary80   & $\varepsilon_e = 2^{-64}$ &  Hardware       \\ 
quadruple & binary128  & $\varepsilon_q = 2^{-113}$ & Software       \\ 
\bottomrule\noalign{\smallskip}
\end{tabular}
  \end{center}
  \caption{The various floating point formats used and their support on common
    hardware. The $\varepsilon$-terms denote the unit roundoff error (for 
    rounding to nearest). We use the letters $s$, $d$, $e$ and
    $q$ synonymously with 32, 64, 80, and 128. For instance,
    $\varepsilon_{_{32}} = \varepsilon_s$.}
  \label{tab:precisions}
\end{table}

In principle, we could perform the entire computation in $y$-bit arithmetic
and, at the end, cast the results to form the $x$-bit output; for 
all practical purposes, we would obtain improved results as
desired. This naive approach, however, is not satisfactory for two reasons:
($i$)~since the eigenvectors need to be
stored explicitly in the
\binaryy\ format, the memory requirement is increased; and more importantly,
($ii$) if the $y$-bit floating point 
arithmetic is much slower than the $x$-bit one, the
performance suffers severely. While the first issue is addressed rather
easily (as discussed Section~\ref{sec:memcost}), the latter requires more
care. The key insight is that it is unnecessary to compute eigenpairs with
residual norms and orthogonality bounded by $\order{n \varepsilon_y}$; instead,
these bounds are relaxed to $\order{\varepsilon_x \sqrt{n}}$ (for
example, think of $\varepsilon_x \approx 10^{-16}$, $\varepsilon_y
\approx 10^{-34}$, and $n \approx 10{,}000$). While in a conventional
implementation the choice of parameters is very restricted, as we show below, we gain enormous freedom in their choice. In
particular, while meeting our new accuracy goals, we are able to select values such
that the amount of necessary computation is reduced, the 
robustness is increased, and parallelism is improved. As our following analysis
shows, we can emphasize the importance of any of those features. 

\subsection{Adjusting the algorithm} Consider the input/output being in a
$x$-bit format and the entire computation being performed in
$y$-bit arithmetic. Starting from this configuration, we expose the new
freedom in the choice of several parameters and justify other changes made to the algorithm.
For example, we identify parts that can be executed in $x$-bit
arithmetic, which might be considerably faster. 

Assuming $\varepsilon_y \ll \varepsilon_x$ (again, think of $\varepsilon_x
\approx 10^{-16}$ and $\varepsilon_y \approx 10^{-34}$), we simplify
Theorem~\ref{resthm} by canceling terms that are insignificant  
even with adjusted parameters (i.e., terms that are comparable to
$\varepsilon_y$ in magnitude\footnote{In particular, we require that $n \varepsilon_y \leq
  \varepsilon_x \sqrt{n}$.}). In our argumentation, we hide all constants,
which anyway correspond to the bounds attainable for a solver purely based
on \binaryy.  For any reasonable implementation of the algorithm, we have the following:
$\alpha = \order{\varepsilon_y}$, $\eta = \order{n \varepsilon_y}$,
$\xi_{\,\downarrow} = \order{\varepsilon_y}$, $\xi_{\uparrow} = \order{\varepsilon_y}$. 
Thus, the orthogonality of the final result is given by 
\begin{equation}
|\hat{z}_i^* \hat{z}_j| = \mathcal{O}\left( k_{rs} \frac{n
    \varepsilon_y}{gaptol} + k_{rr} \, d_{max}\, \frac{n \varepsilon_y}{gaptol} \right) \,.
\label{simpleorthobound}
\end{equation}
Similarly, for the bound on the residual norm, we get 
\begin{equation}
\norm{M_{root}\, \hat{z}_i - \hat{\lambda}_i[M_{root}]  \, \hat{z}_i} =
\mathcal{O}\left(\norm{r^{(local)}} + \gamma \, spdiam[M_{root}]\right)
\label{simpleresbound}
\end{equation}
with $\norm{r^{(local)}} \leq k_{rs} \, gap\left(\hat{\lambda}_i[M]\right) \frac{n
\varepsilon_y}{gaptol}$
 and $\gamma = \order{k_{elg} \, d_{max} \, n \varepsilon_y}$.

We now provide a list of changes that can be done to the algorithm. 
We discuss their effects on performance, parallelism, and memory
requirement.  

\paragraph{Preprocessing}
We assume scaling and splitting is done as in a solver purely based on $x$-bit
floating point arithmetic. In particular, off-diagonal
element $e_i$ of the input, $1 \leq i \leq n-1$, is set to zero whenever
\begin{equation*} 
|e_i| \leq \varepsilon_x \norm{T} \,,
\end{equation*}  
where $n$ and $T$ refer to the {\it unreduced} input.\footnote{We can
  relax the condition further and use
  $|e_i| \leq \varepsilon_x \sqrt{n} \norm{T}$.} 
We remark that this
criterion is less strict than setting elements to zero whenever $|e_i| \leq
\varepsilon_y \norm{T}$. Splitting the input matrix into submatrices is
beneficial for both performance and accuracy as these are mainly determined
by the largest submatrix.
In the rest of this
section, we assume that the preprocessing has been done and each subproblem
is treated independently by invoking Algorithm~\ref{alg:mrrr}. In particular, whenever
we refer to matrix $T$, it is assumed to be irreducible; whenever we
reference the matrix size $n$ in the context of parameter settings, it refers
to the size of the processed block. 

\paragraph{Choice of representations} For different forms
of representing tridiagonals (e.g., bidiagonal, twisted, or blocked factorizations) and their
data (e.g., $N$-, $e$-, or $Z$-representation), different algorithms
implement the shift operation: $M_{shifted} = M - \tau I$. All these
algorithms are stable in the sense that the relation holds 
exactly if the data for $M_{shifted}$ and $M$ are perturbed element-wise by a
relative amount bounded by $\order{\varepsilon_y}$. The implied constants
for the perturbation bounds vary slightly. As $\varepsilon_y \ll
\varepsilon_x$, instead of concentrating on accuracy issues, we make our
choice based on robustness and {\it performance}. A discussion of performance
issues related to different forms of the representations can be found
in~\cite{Willems:twisted,Willems:Diss}. Based on this discussion, it appears
that twisted factorizations with $e$-representation seem to be a reasonable
choice. As the off-diagonal entries of all the matrices stay the same, they
only need to be stored once and are reused during the entire
computation. 

\paragraph{Random perturbations} In Line~\ref{line:mrrr:perturb} of
Algorithm~\ref{alg:mrrr}, to break up tight clusters, the data of 
$M_{root}$, $\{x_1,...,x_{2n-1}\}$, is perturbed element-wise by small
random relative amounts:\footnote{True randomness is not necessary; any
  (fixed) sequence of pseudo-random numbers can be used.} 
$\tilde{x}_i = x_i (1 + \xi_i)$ with $|\xi_i| \leq \xi$ for all $1 \leq i
\leq {2n-1}$. In practice, a value like $\xi = 8
\varepsilon$ is used. Although our data is in \binaryy, we are quite
aggressive and adopt $\xi = \varepsilon_x$ or a small multiple
of it. Thus, for $y = 2x$, about half of the digits in each entry of the representation are chosen
randomly; therefore, with high probability, eigenvalues do not agree to
many more than $-\lceil\log_{10} \varepsilon_x\rceil$ digits.
This has two major effects: ($i$) together with the changes in
$gaptol$ (see below), in practice, the probability to encounter $d_{max} > 1$ 
becomes very low, and ($ii$) it becomes easier to find suitable shifts such that the
resulting representation satisfies the requirements of relative robustness
and conditional element growth. The positive impact of small 
$d_{max}$ on the accuracy is apparent from \eqref{simpleorthobound}
and \eqref{simpleresbound}. Furthermore, as 
discussed below, due to limiting $d_{max}$, the computation can be
reorganized for efficiency. Although it might look innocent, the more
aggressive random perturbations lead to much improved robustness: A detailed
discussion can be found in~\cite{Dhillon05gluedmatrices}.\footnote{For a quantitative
  assessment of robustness, see~\cite{mydiss}.} 

\paragraph{Classification of the eigenvalues} Due to the importance of the
$gaptol$-parameter, adjusting it to our requirements is key to the success
of our approach. The parameter influences nearly all stages of the
algorithm; most importantly, the classification of eigenvalues into
well-separated and clustered. As already discussed, the choice of {\it
  gaptol} is restricted by the loss of orthogonality that we are willing to
accept; in practice, the value is often
chosen to be $10^{-3}$~\cite{Dhillon:2004:MRRR}.\footnote{For instance,
  LAPACK's {\tt DSTEMR} uses $10^{-3}$, while {\tt SSTEMR} uses $3 \cdot
  10^{-3}$.} As we merely require 
orthogonality of $\varepsilon_x \sqrt{n}$, we accept
more than three orders of magnitude loss of orthogonality. Both terms in
\eqref{simpleorthobound} (and the in practice observed orthogonality)
are proportional to $n 
\varepsilon_y / gaptol$. This means that the value of {\it gaptol} can be chosen as small as $
 \varepsilon_y \sqrt{n} / \varepsilon_x$. As a consequence, we might select any
 value satisfying
\begin{equation}
 \min\left\{ 10^{-3}, \frac{\varepsilon_y \sqrt{n}}{\varepsilon_x} \right\} \leq gaptol \leq
 10^{-3} \,,
\label{eq:gaptolinterval}
\end{equation}
 where the $10^{-3}$ terms are derived from practice and might be altered slightly.
 Note that $gaptol$ potentially becomes as small
 as $10^{-9} \sqrt{n}$ in the single/double case and $10^{-18} \sqrt{n}$ in the
 double/quadruple one. If we restrict the analysis to matrices with size $n \leq 10^6$, we can
 choose a constant $gaptol$ as small as $10^{-6}$ and $10^{-15}$ 
 respectively for the single/double and double/quadruple cases. 
 
 With any choice of $gaptol$ complying \eqref{eq:gaptolinterval}, accuracy
 to the desired level is guaranteed, and there is room to choose the specific
 value of $gaptol$, as well as other parameters, to optimize performance or
 parallelism. In particular, by generally reducing the clustering of the
 eigenvalues, the smallest possible value of $gaptol$ provides the greatest parallelism.
 To quantify this statement, for any matrix, we
 define {\it clustering} $\rho \in [1/n, 1]$ formally as the size of the
 largest cluster divided by the matrix size. There are two main advantages
 in decreasing $\rho$: ($i$) the work is reduced as processing the
 largest cluster introduces $\order{\rho n^2}$ 
 flops extra work, and ($ii$) the potential parallelism is increased. A conservative
 estimate of the parallelism of a problem is provided by $\rho^{-1}$. For instance, $\rho
 = 1/n$ implies that the problem is {\it embarrassingly parallel}. The
 estimate of parallelism assumes that clusters are processed sequentially,
 while in reality the bulk of the work (the refinement of the eigenvalues and
 the final computation of eigenpairs) can be
 parallelized. Nonetheless, matrices with high clustering still pose
 difficulties to MRRR as they introduce load-balancing issues and
 communication, which considerably reduce the parallel
 scalability~\cite{Vomel:2010:ScaLAPACKsMRRR,VoemelRefinedTree2007tr,mydiss}. 
 Therefore, even if we did not have the desire to guarantee improved accuracy of
 the method, we could use the mixed precision approach to significantly
 enhance parallelism. In this case, the $\sqrt{n}$-dependence on the lower
 bound for the value of $gaptol$ would be removed and the bound could be
 loosened by another three orders of magnitude; that is, we could choose a value of
 $10^{-12}$ and $10^{-21}$ for the single/double and double/quadruple case,
 respectively.\footnote{If we select values $10^{-9}$ and $10^{-18}$, we improve the bounds by three orders of magnitude.} Consequently, {\it almost
 all} computations become embarrassingly parallel. 

 As an example, Table~\ref{tab:clustering}
 shows the clustering for double precision Hermite 
 type\footnote{See~\cite{Marques:2008} for information on test matrices.} test matrices of
 various sizes with four distinct classification criteria:\footnote{Criterion I is used in LAPACK~\cite{DesignMRRR} and in results of {\tt
    mr3smp} in~\cite{mr3smp}, which usually uses II. Criterion II is used in
  ScaLAPACK~\cite{Vomel:2010:ScaLAPACKsMRRR} and
  Elemental~\cite{EleMRRR}. In massively parallel computing environments,
  criteria III and IV can (and should) additionally complemented with the
  splitting based on absolute gaps; see also~\cite{mixedtr}.} 
(I) $gaptol = 10^{-3}$, (II) $gaptol = 10^{-3}$, combined with splitting based on
the absolute gap as proposed in~\cite{VoemelRefinedTree2007tr} to enhance
parallelism, (III) $gaptol = 10^{-10}$, and (IV) $gaptol = 10^{-15}$.
\begin{table}[htb]
\begin{center}
\small
\begin{tabular}{c@{\quad\quad\quad}cccc}
\toprule
Criterion &   \multicolumn{4}{c}{Matrix size} \\
               & 2{,}500 & 5{,}000 & 10{,}000 & 20{,}000 \\
\midrule
I  &  0.70  &         0.86   &  0.93  & 0.97 \\
II &  0.57   &        0.73   &  0.73  & 0.73 \\
III &  4.00e-4  &   2.00e-4   &  1.00e-4  & 5.00e-5  \\
IV &  4.00e-4  &   2.00e-4   &  1.00e-4  & 5.00e-5 \\
% I  &  1746  & 4292 & 9309 & 19316 \\
% II &  1434  & 3632 & 7347 & 14691 \\
% III &  1        &   1     &  1      & 1  \\
% IV &  1        &   1     &  1      & 1 \\
\bottomrule\noalign{\smallskip}
\end{tabular} 
\end{center}
\caption{The $gaptol$-parameter effect on clustering $\rho \in [1/n,1]$.  
}
\label{tab:clustering}
\end{table}
For the latter two criteria, the computations are embarrassingly parallel.
As with this example, experience shows that, thanks to a reduced value of
$gaptol$ as in criteria III or IV, many problems 
become embarrassingly parallel {\em and} guarantee improved accuracy. 
In case $\rho = 1/n$, $d_{max}= 0$, which not only benefits accuracy by
\eqref{simpleorthobound} and \eqref{simpleresbound}, but also has a 
more dramatic effect: {\it the danger of 
not finding representations that satisfy the requirements is entirely
removed.} This follows from the fact that a satisfactory root representation is 
always found (e.g., by making $T - \mu I$ definite) and no other
representation needs to be computed. 
Even in cases with $d_{max} > 0$, the
number of times Line~\ref{line:mrrr:shifting} of Algorithm~\ref{alg:mrrr}
needs to be executed is often considerably reduced. 

On the downside, selecting a smaller $gaptol$ can result in more work in the
initial approximation\footnote{For instance, if bisection is used to obtain
  initial approximations to the eigenvalues.} and later refinements -- in
both cases, eigenvalues must be approximated to relative accuracy of about $gaptol$, see
\eqref{eq:relaccgaptol}; hence, optimal performance is often not achieved
for the smallest possible value of $gaptol$. 
Moreover, as we discuss below, if one is willing to limit
the choice of $gaptol$, the computation and refinement of eigenvalues can be
done (almost) entirely in $x$-bit arithmetic.\footnote{For the refinement of extreme
  eigenvalues prior to selecting shifts, we still need to resort to
  $y$-bit arithmetic.} If $y$-bit arithmetic is slow, it might be best to
take advantage of the faster $x$-bit arithmetic. And, as we see below as
well, if not the smallest possible value is chosen for $gaptol$, the
requirements the intermediate 
representations must fulfill are relaxed, thereby increasing the robustness
of the method.  

Another corollary of adjusting $gaptol$ is slightly hidden: in
Line~\ref{line:mrrr:shifting} of Algorithm~\ref{alg:mrrr}, we gain more freedom
in selecting $\tau$ such
that, at the next iteration, the index set 
$\mathcal{I}_r$ splits into two or more subsets. For instance, when choosing
$\tau$ close to one end of the cluster, we are able to ``back off'' further
away than usual from the end of the cluster in cases where, in a previous
attempt, we did not find a representation satisfying the requirements~\cite{DesignMRRR}. 

We cannot overemphasize the positive effects an adjusted {\em gaptol} has on
robustness and parallel scalability. In particular, in a massively parallel
computing environment, the smallest value for $gaptol$ significantly
improves the parallel scalability. And since many problems become embarrassingly
parallel, the danger of failing to find suitable representations is entirely
removed.

\paragraph{Arithmetic used to approximate eigenvalues} 
In Lines~\ref{line:mrrr:initialeigvals} and \ref{line:mrrr:refine} of
Algorithm~\ref{alg:mrrr}, eigenvalues are respectively computed and refined
to a specified relative accuracy. 
In both cases, we are given a
representation, which we call $M_y$ henceforth, and an index set
$\mathcal{I}$ that indicates the eigenvalues that need to be approximated. 
When the $y$-bit arithmetic is much slower than the $x$-bit one (say a
factor 10 or more), the use of the latter is preferred: 
One creates a temporary copy of $M_y$ in \binaryx\ -- called $M_x$ henceforth
-- that is used for the eigenvalue computation in $x$-bit arithmetic. The creation of $M_x$
corresponds to an element-wise relative perturbation of $M_y$ bounded by
$\varepsilon_x$. By the relative robustness of the representation, 
\begin{equation}
\abs{\lambda_i[M_x] - \lambda_i[M_y]} \leq
k_{rr} n \varepsilon_x \abs{\lambda_i[M_y]} \,.
\end{equation}
For instance, bisection can be used to compute eigenvalue approximations
$\hat{\lambda}_i[M_x]$ to high relative accuracy, after which $M_x$ is discarded. As
casting the result back to \binaryy\ causes no additional error, it is $\hat{\lambda}_i[M_y] =
\hat{\lambda}_i[M_x]$ and
\begin{equation*}
\abs{\hat{\lambda}_i[M_y] - \lambda_i[M_x]} \leq
k_{bi} n \varepsilon_x \abs{\lambda_i[M_x]} \,,
\end{equation*}
where $k_{bi}$ is a moderate constant given by the bisection method. To first order, by the
triangle inequality, it holds
\begin{equation}
\abs{\hat{\lambda}_i[M_y] - \lambda_i[M_y]} \leq \left( k_{rr} + k_{bi}\right) n
\varepsilon_x \abs{\lambda_i[M_y]} \,.
\label{eq:accuracyeigvalinxarithmetic}
\end{equation}
Provided $\left( k_{rr} + k_{bi}\right) n \varepsilon_x \lesssim
gaptol$, by \eqref{eq:relaccgaptol}, $x$-bit arithmetic can be used to approximate the
eigenvalues. Thus, an additional constraint on both the size $n$ and
$gaptol$ arises:
Given a $gaptol$, we must limit the matrix size up to which we do the
computation purely in $x$-bit arithmetic. Similarly, for a given matrix size, we need to adjust
the lower bound on $gaptol$ in \eqref{eq:gaptolinterval}. As an example, if
say $k_{rr} \leq 10$, $k_{bi} \leq 10$, $n \leq 10^6$, and $\varepsilon_x =
\varepsilon_d = 2^{-53}$, it is required that that $gaptol \gtrsim 10^{-10}$. 
When resorting to $x$-bit arithmetic or if $gaptol$ is chosen too small, one might
respectively verify or refine the result of the $x$-bit eigenvalue computation using
$y$-bit arithmetic without significant costs.\footnote{If the first
  requirement in Definition~\ref{def:RRR} is 
  removed, we can still make use of $x$-bit arithmetic although
  \eqref{eq:accuracyeigvalinxarithmetic} might
  not always be satisfied anymore.} 

\paragraph{Requirements on the representations}
As long as $k_{elg} n\varepsilon_y \ll \varepsilon_x \sqrt{n}$, by \eqref{simpleresbound},
the residual with respect the $M_{root}$ is mainly influenced by the local
residual. In our mixed precision approach, without loss of accuracy, it
is possible to allow for  
\begin{equation}
k_{elg}\leq \max \left\{10, \frac{\varepsilon_x}{\varepsilon_y \sqrt{n}}
\right\} \,,
\label{eq:newkelgbound}
\end{equation}
where we assumed 10 was the original value of $k_{elg}$. As a result, the
requirement on the conditional element growth is considerably relaxed. 
For instance, in the single/double and double/quadruple cases,
assuming $n \leq 10^6$, bounds on $k_{elg}$ of about $10^6$ and $10^{15}$ are
sufficient, respectively. 
If $gaptol$ is
not chosen as small as possible, the bound on $k_{rr}$ is loosened in a
similar fashion:
\begin{equation}
k_{rr}\leq \max \left\{10, \frac{\varepsilon_x}{\varepsilon_y
    \sqrt{n}} \cdot gaptol \right\} \,.
\label{eq:newkrrbound}
\end{equation}
As an example, in the double/quadruple case, assuming $n \leq 10^6$ and
$gaptol$ set to $10^{-10}$, $k_{rr} \leq 10^5$ would be sufficient to ensure
accuracy. 

\paragraph{Rayleigh quotient iteration}
Our willingness to lose orthogonality up to a certain level, which is
noticeable in the lower bound on $gaptol$, 
is also reflected in \eqref{localresbound}. 
As $n \varepsilon_y / gaptol \leq \varepsilon_x 
 \sqrt{n}$, we stop the RQI when
\begin{equation}
\norm{r^{(local)}} \leq k_{rs}  \cdot
gap\left(\hat{\lambda}_i[M]\right) \, \varepsilon_x \sqrt{n} \,,
\end{equation}
where $k_{rs}$ is $\order{1}$. In practice, we take $k_{rs} = 1$ (or
even $k_{rs} = 1/\sqrt{n}$). As a consequence, the iteration is stopped
earlier on, thereby reducing the overall work. 

As a side note: In the rare cases where RQI fails to converge (or as a general alternative to
RQI), we commonly resort to bisection to approximate the eigenvalue
$\lambda_i$ and then use only one step of RQI (with or without
applying the correction term). In the worst case, we require the eigenvalue to be
approximated to high relative accuracy, $|\hat{\lambda}_i -
\lambda_i| = \order{n \varepsilon_y
  |\lambda_i|}$~\cite{Dhillon:2004:Ortvecs}. With mixed precision, we 
relax the condition to $|\hat{\lambda}_i - \lambda_i| =
\order{\varepsilon_x \sqrt{n} |\lambda_i| \, gaptol}$, which is less
restrictive if $gaptol$ is not chosen as small as
possible.\footnote{The implied constants being the same and given by the
  requirement of a regular solver based on $y$-bit arithmetic. In a similar way, we
  could say that the Rayleigh quotient correction does not
  improve the eigenvalue essentially anymore 
  if $|\gamma_s|/\norm{\hat{z}_i} = \order{\varepsilon_x |\hat{\lambda}_i|
    \, gaptol / \sqrt{n}}$, instead of $|\gamma_s|/\norm{\hat{z}_i} =
  \order{\varepsilon_y |\hat{\lambda}_i|}$. We never employed it as such a
  change will hardly have any effect on the computation time.} If
$relgap(\hat{\lambda}_i) \gg gaptol$, the restriction on the accuracy of the
approximated eigenvalue is lifted even further~\cite{Willems:Diss}. 

\paragraph{The representation tree} 
Thanks to the random
perturbation of the root representation and a properly adjusted
$gaptol$-parameter, we rarely expect to see large values for $d_{max}$. For all
practical purposes, in the case of $y = 2x$, we may assume $d_{max} \leq
2$. As a result,
the computation can be rearranged, as discussed
in~\cite{Willems:framework} and summarized in the following: To bound the memory
consumption, a breath-first strategy such as in Algorithm~\ref{alg:mrrr} is used; see for
instance in~\cite{DesignMRRR,mr3smp}. This means that, at
any level of the representation tree, all singletons are processed before the
clusters. A depth-first strategy would instead process entire clusters,
with the
only disadvantage that meanwhile up to $d_{max}$ representations need to 
be kept in memory. If $d_{max}$ is limited as in our case, the depth-first
strategy can be used without disadvantage. In fact, a depth-first strategy
brings two advantages: ($i$) 
copying representations to and from the eigenvector matrix is avoided entirely (see
Section~\ref{sec:memcost} on the benefit for the mixed precision approach) and ($ii$) if no
suitable representation is found, there is the possibility of
backtracking, that is, we process the cluster again by choosing different
shifts at a higher level of the representation tree. For these reasons, in the mixed
precision approach, a depth-first strategy is preferred. 

\subsection{Memory cost}
\label{sec:memcost}

We stress that in our approach, both input and output are in \binaryx\
format; only {\it internally} (i.e., hidden to a user) $y$-bit
arithmetic is used.
The memory management of an actual implementation of MRRR is affected by the
fact that matrix $Z \in \mathbb{R}^{n \times k}$, which on output contains
the desired eigenvectors, 
is commonly used as intermediate work space. Since $Z$ is in \binaryx\ format, whenever $y >
x$, the work space is not sufficient anymore for its customary use: For
each index set $\mathcal{I}_r$ with $|\mathcal{I}_r| > 1$, a
representation, $M_{shifted}$, is stored in the corresponding columns of
$Z$~\cite{DesignMRRR,mr3smp}. 
As these representations consist of $2n-1$ \binaryy\ numbers, this approach
is generally not applicable. However, if we restrict to $y \leq 2x$, we can store the 2$n$ 
\binaryy\ numbers whenever a cluster of size four and more is encountered. %$|\mathcal{I}_r| \geq 4$. 
Thus, the computation
must be reorganized so that at least clusters containing less than four
eigenvalues are processed without storing any data in $Z$ temporarily. In
fact, using a depth-first strategy, we remove the need to use $Z$ as
temporary storage entirely. 
Immediately after computing an eigenvector
in \binaryy, it is converted to \binaryx, written into $Z$, and
discarded. While our approach slightly increases the memory usage, we
do not require much more memory: with $p$ denoting the number of computational
threads, our mixed precision solver still needs only $\mathcal{O}(p n)$
\binaryx\ floating point numbers extra work space.

\section{Practical aspects}
\label{sec:implementation}

We have implemented the mixed precision approach for three cases:
{\it single/double}, {\it double/extended}, and {\it double/quadruple}. The
first solver accepts single precision input and produces single precision
output, but internally uses (hidden to the user) double precision. The
other two are for double precision input/output. The performance of the
solvers, compared with the traditional implementation, depends entirely on
the difference in speed between the two involved arithmetic. 
If the higher precision arithmetic is not much slower (say less
than a factor four), the approach is expected to always work well, even for
sequential executions and relatively small matrices. If the higher
precision arithmetic is considerably slower, the mixed precision approach
might still perform well for large 
matrices. Due to increased parallelism, our approach is also
expected to perform generally well on highly parallel systems. 
Our target application is the computation of a subset of eigenpairs of
large-scale dense Hermitian matrices. For such a scenario, we 
tolerate a slowdown of the tridiagonal eigensolver due to mixed precisions
without affecting performance significantly as the reduction to tridiagonal
form is the performance bottleneck~\cite{EleMRRR,mixedtr}.

\subsection{Implementations}
In Section~\ref{sec:experiments}, we present experimental results of our
implementations. All mixed precision implementations are based on a 
multi-threaded variant of MRRR, {\tt mr3smp}, presented
in~\cite{mr3smp,para2010}, which is built on top of LAPACK's routine {\tt
  DSTEMR} (version 3.2). All codes use
$N$-representations of lower bidiagonal 
factorizations. Bisection is used
for the initial eigenvalue computation if a small 
subset of $k$ eigenpairs is requested or if the
number of executing threads exceeds $12k/n$~\cite{DesignMRRR,mr3smp}. If
all eigenpairs are requested and the number of threads is less than 12, the fast
sequential {\it dqds algorithm}~\cite{AccurateSVDandQDtrans,dqds99} is used
instead of bisection. As a consequence, speedups compared with the sequential
execution appear less than perfect even for an embarrassingly parallel
computation. 

We did not relax the requirements on the representations
according to \eqref{eq:newkelgbound} and \eqref{eq:newkrrbound}; we only
benefit from the possibility of 
doing so indirectly: If no suitable 
representation is found, the best candidate is chosen, which might still
fulfill the relaxed requirements.  

In the following, we provide
additional comments to all of the mixed precision solvers individually. 
As parameters can take a wide range of values (in particular, $gaptol$, but also $k_{rr}$ and
$k_{elg}$) and several design decisions can be made, optimizing a code for performance
is non-trivial as it generally depends on both the specific input and the architecture. 
While we cannot expect to create an ``optimal''
design for all input matrices and architectures, we make design
decisions in a way that in general yields good performance. 
For instance, on
a highly parallel machine one would select a small value for $gaptol$ to
increase parallelism. For testing purposes, we disabled the
classification criterion based on the absolute gaps of the eigenvalues
proposed in~\cite{VoemelRefinedTree2007tr}, which might reduce clustering even
further (it has no consequences for our test cases shown in the next section). 

\paragraph{Single/double} With widespread language and hardware support for
  double precision, the mixed precision approach is most easily implemented
  for the {\it single/double} case. In our test implementation, we fixed
  $gaptol$ to $10^{-5}$. When bisection is used, the initial eigenvalue
  approximation is done to a relative accuracy of $10^{-2} \cdot gaptol$. As
  on most machines the double precision arithmetic is not more than a factor two
  slower than the single precision one, we carry out {\it all} computations in
  the former. Data conversion is only necessary when reading the input and
  writing the output. As a result, compared with a double precision solver
  using a depth-first strategy,
  merely a number of convergence criteria and thresholds must be
  adjusted, and the RQI must be performed using a temporary vector that is,
  after convergence, written into the output eigenvector matrix. The mixed precision code
  closely resembles a conventional double precision implementation of MRRR.
\paragraph{Double/extended} Many current architectures have hardware
  support for a 80-bit extended floating point format (see
  Table~\ref{tab:precisions}). As the unit roundoff $\varepsilon_e$ is only
  about three orders of magnitude smaller than $\varepsilon_d$, we can 
  improve the accuracy of MRRR by this amount. For matrices of moderate size, this means
  that the accuracy becomes comparable to that of the best methods. The main
  advantage of the extended format is that, compared with double precision, its
  arithmetic comes without any or only a small loss in speed. However, we
  cannot make any further adjustments in the algorithm, which 
  positively effect its robustness and parallelism. 
  We do not include results for the {\it double/extended}
case in the next section; however, we tested the approach and experimental
results can be found in~\cite{mydiss,mixedtr}.  
\paragraph{Double/quadruple} As quadruple precision arithmetic is not
  widely supported by today's processors or languages, we had to resort to
  software-simulated arithmetic, which is rather slow. For this reason, we
  used double precision for the initial approximation and for the refinement
  of the eigenvalues. The necessary intermediate data conversions make the
  mixed precision approach slightly more complicated to implement than
  the {\it single/double} one. 
  We used the value $10^{-10}$ for $gaptol$ in our tests. Further details
  can be found in~\cite{mixedtr}.  

\subsection{Portability} The biggest problem of the mixed precision approach
is a potential lack of support for the involved data types. As single and
double precisions are supported by 
virtually all machines, languages, and compilers, the mixed precision
approach can be incorporated to any linear algebra library for single
precision input/output. However, for double precision input/output, we
need to resort to either 
extended or quadruple precision.  
Not all architectures, languages, and compiler support these formats. For
instance, the 80-bit floating point format is not supported by all
processors. Futhermore, while the FORTRAN {\tt REAL*10} data
type is a non-standard feature of the language and is not supported by all
compilers, a C/C++ code can use the standardized {\tt long double} data type
(introduced in ISO C99) that achieves the 
desired result on most architectures that support 80-bit arithmetic. 
For the use of quadruple precision, there are presently two major drawbacks:
($i$) it is usually not supported in hardware, which means that one has to resort to
a rather slow software-simulated arithmetic, and ($ii$) the support from compilers  
and languages is rather limited. 
While FORTRAN has a {\tt REAL*16}
data type, the quadruple precision data type in C/C++ is compiler-dependent:
  for instance, there exist the {\tt \rule{8pt}{0.5pt}float128} and {\tt
  \rule{4pt}{0.5pt}Quad} data types for the GNU and Intel
compilers, respectively. An external library implementing the
software arithmetic might be used for portability. In all cases, the
performance of quadruple arithmetic depends on its specific 
implementation. It is however likely that the hardware/software support for
quadruple precision will be improved in the near future.

\section{Experimental Results}
\label{sec:experiments}

All tests, in this section, were run on a multi-processors system comprising four
eight-core {\it Intel Xeon X7550 Beckton} processors, with a nominal clock
speed of 2.0 GHz. 
Subsequently, we refer to this
machine as {\sc Beckton}. We used LAPACK version
3.4.2 and linked the library with the vendor-tuned MKL BLAS version
12.1. In addition to the results for LAPACK's routines and our mixed
precision solvers, we also include results for {\tt
mr3smp}~\cite{mr3smp}. All routines were compiled with
Intel's compiler version 12.1 and optimization level {\tt -O3} enabled.
Although we present only results for computing {\it all}
eigenpairs (LAPACK's DC does not allow the computation of subsets), we
mention that {\it MRRR's strength and main application lies in the
  computation of subsets of 
eigenpairs}. 

For our tests, we used matrices of size ranging from $2{,}500$ to
$20{,}000$ (in steps of $2{,}500$) of six different types: uniform eigenvalue distribution,
geometric eigenvalue distribution, 1--2--1, Clement, Wilkinson, and Hermite. The
dimension of the Wilkinson type matrices is $n+1$, as they are
only defined for odd sizes. Details on these matrix types can be found
in~\cite{Marques:2008}. 
To help the exposition of the results, in the accuracy plots, the matrices
are sorted by type first and then by size; vice versa, in the plots relative
to timings, the matrices are sorted by size first and then by type.  

\begin{figure}[bth]
   \centering   
   \subfigure[Execution time: sequential.]{
     \includegraphics[width=.47\textwidth]{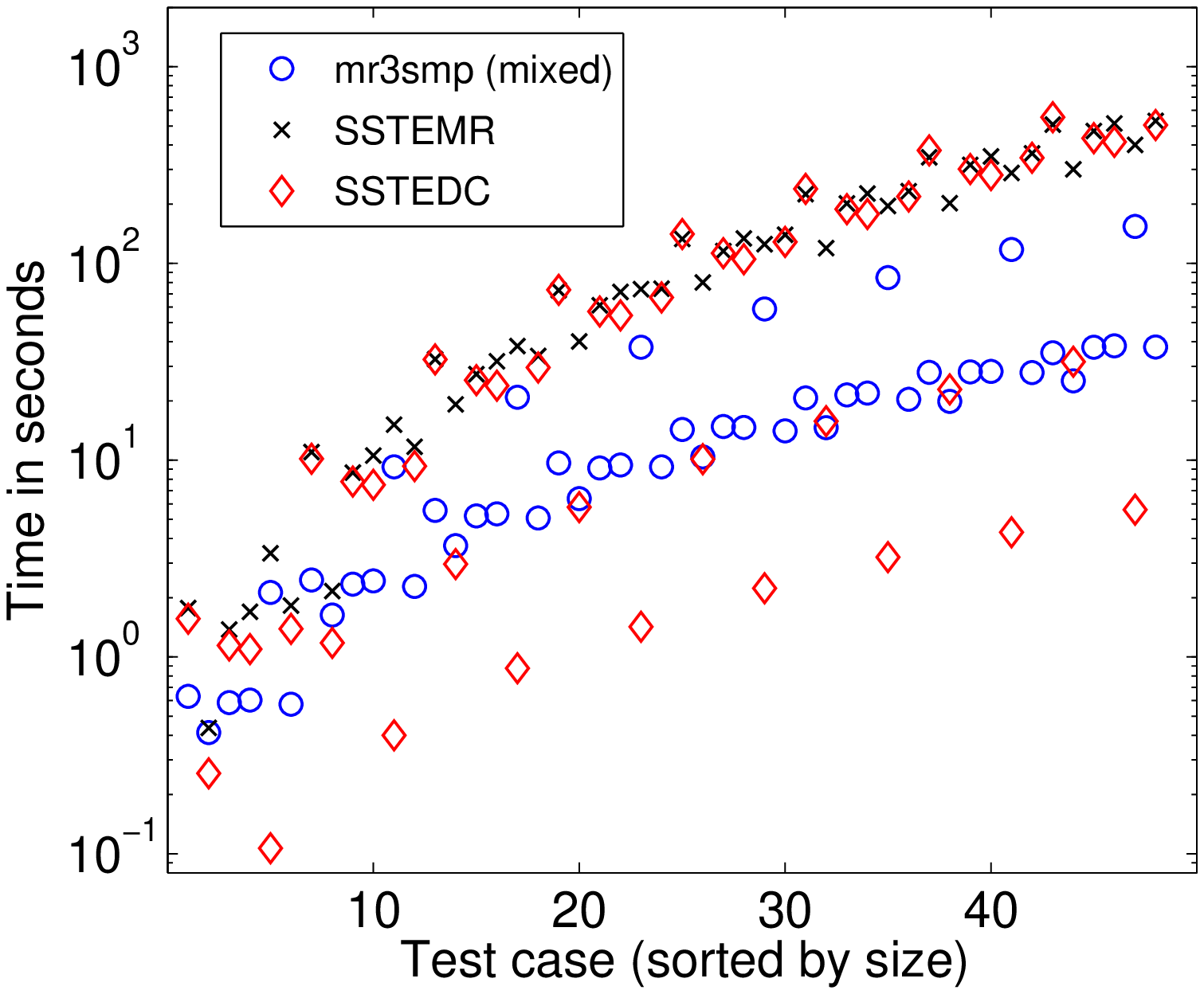}
   \label{fig:accartifial32coresinglea}
   } \subfigure[Execution time: multi-threaded.]{
     \includegraphics[width=.47\textwidth]{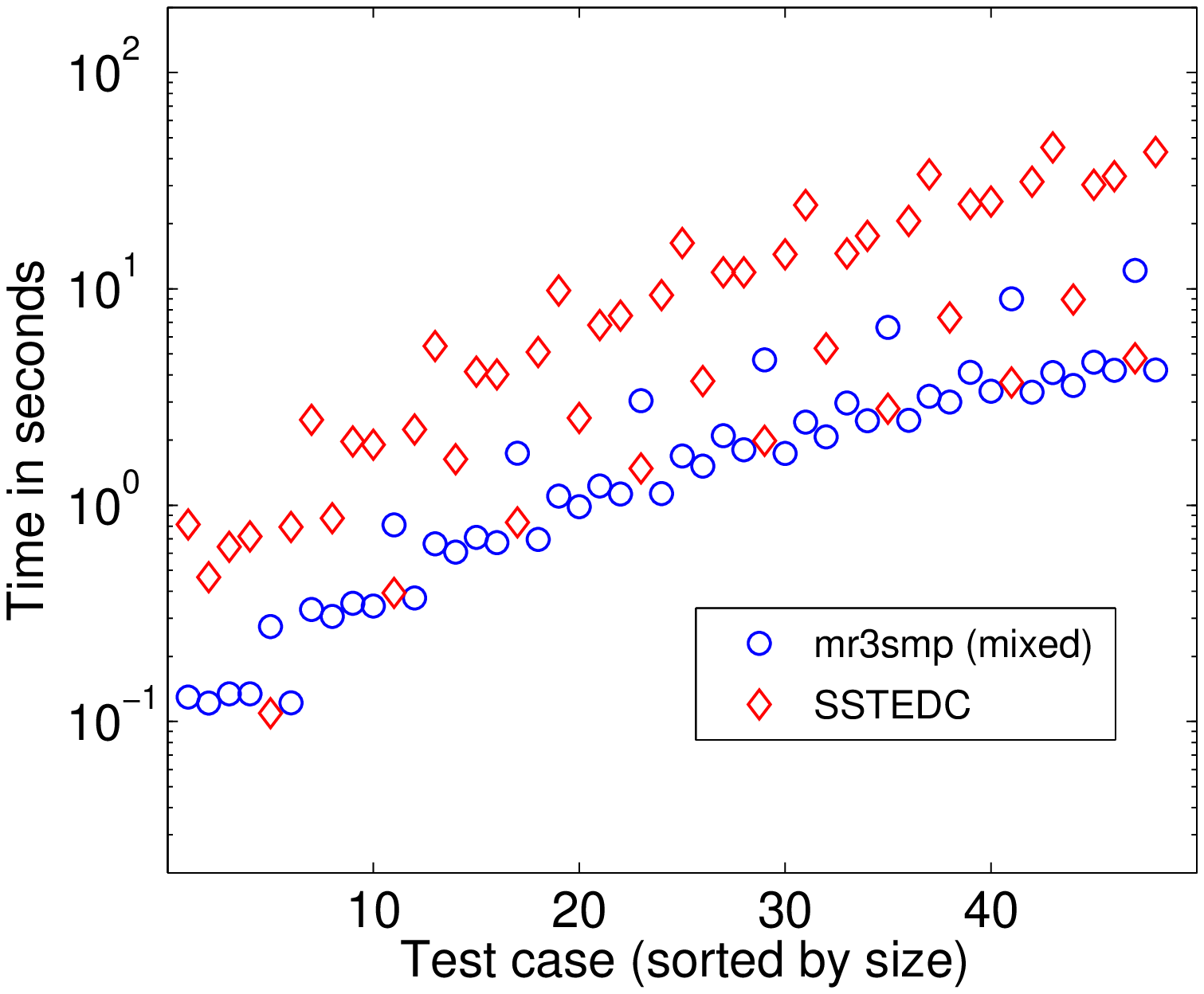}
   \label{fig:accartifial32coresingleb}
   }

   \subfigure[Largest residual norm.]{
     \includegraphics[width=.47\textwidth]{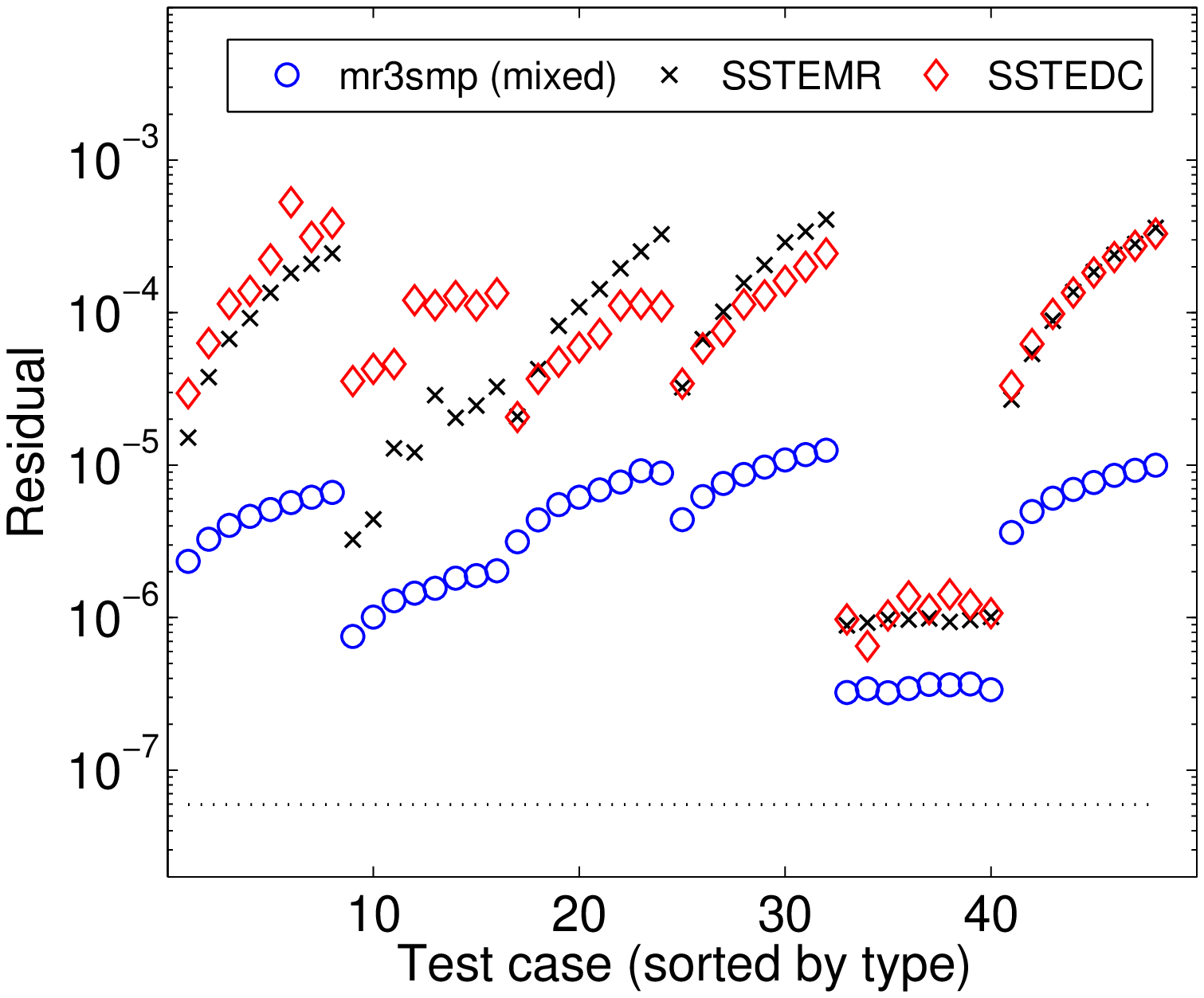}
   \label{fig:accartifial32coresinglec}
   } \subfigure[Orthogonality.]{
     \includegraphics[width=.47\textwidth]{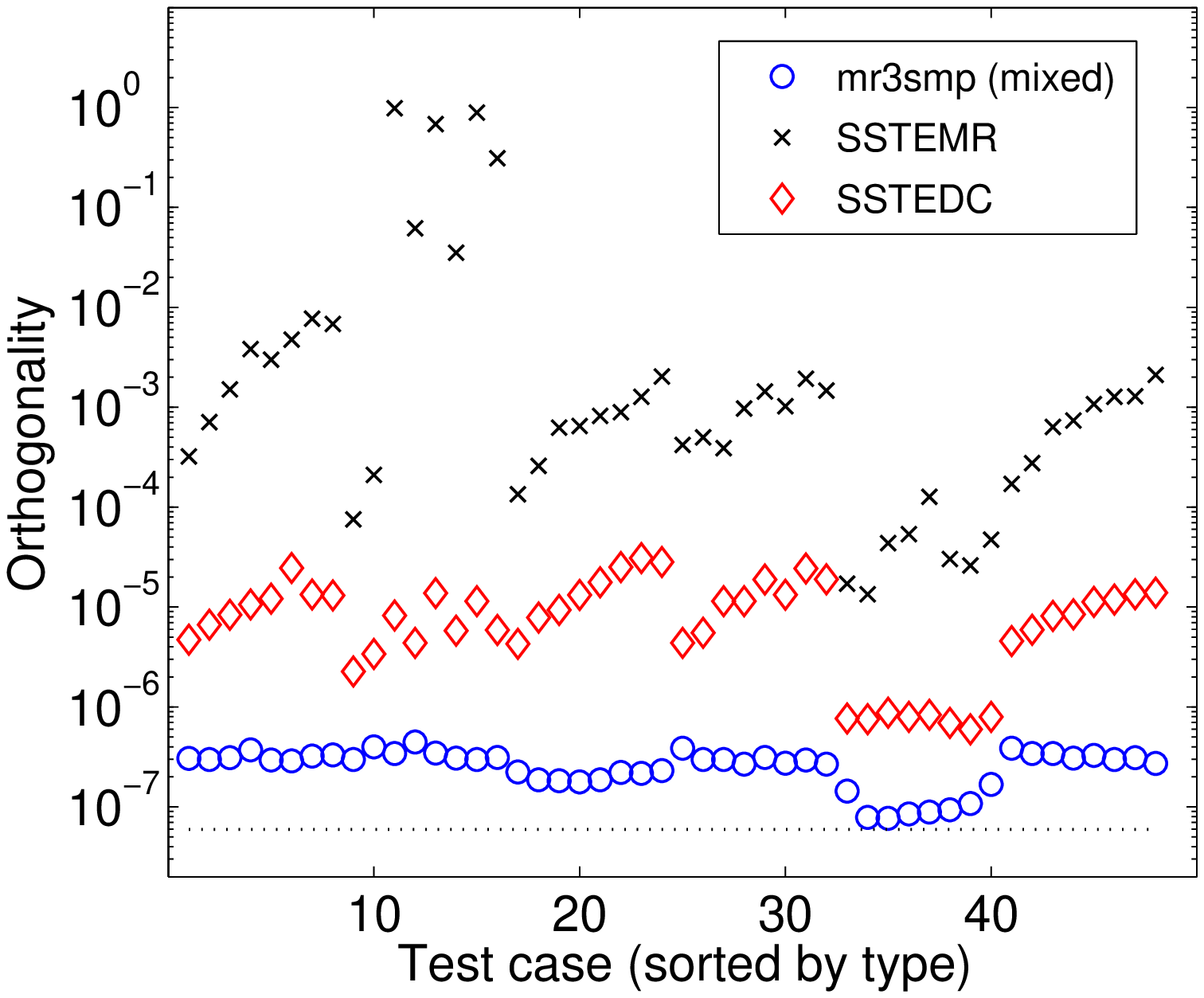}
   \label{fig:accartifial32coresingled}
   }

   \caption{
     Time and accuracy on {\sc Beckton}. Timings are presented in a
     logarithmic scale. The largest residual norm and the orthogonality are
     measured as in~\eqref{def:defresortho}. The dotted black line
     corresponds to unit round-off $\varepsilon_s$. As there exist no parallel MRRR for
     single precision, we show timings for our mixed precision approach and
     {\tt SSTEDC} only. 
   }
   \label{fig:accartifial32coresingle}
\end{figure}

Figure~\ref{fig:accartifial32coresingle} shows
timings and accuracy for single precision inputs. 
As a reference, we include
results for LAPACK's {\tt SSTEMR} (MRRR) and {\tt SSTEDC} (Divide \& Conquer). 
As shown in Fig.~\ref{fig:accartifial32coresinglea}, even in a sequential execution, our mixed precision approach is up to an
order of magnitude {\it faster} than LAPACK's {\tt SSTEMR}. For one type of
matrices, {\tt SSTEDC} 
is considerably faster than for all 
the others. These are the Wilkinson matrices, which represent a class of
matrices that allow for heavy deflation within the Divide \& Conquer
approach. For all other matrices, which do 
not allow such extensive deflation, our solver is {\it faster} than {\tt
  SSTEDC}. As seen in Fig.~\ref{fig:accartifial32coresingleb}, in a
parallel execution, the performance gap for the Wilkinson matrices almost
entirely vanishes, while for the other matrices our solver remains up to an order of magnitude
faster than {\tt SSTEDC}.  
As depicted in Figs.~\ref{fig:accartifial32coresinglec}--\subref{fig:accartifial32coresingled},
our routine is not only as accurate as
desired but it is the most accurate one. For single precision input/output
arguments, we obtain a solver that is more accurate {\it and} faster than the
original single precision solver. In addition,
the solver is more scalable, and more robust. In 38 out of the 48 test cases,
{\tt SSTEMR} accepted representations that did {\em not} pass the test for
relative robustness, thereby jeopardizing the accuracy of the
result. In contrast, using mixed precisions, our solver was able to find
suitable representations in all cases.

\begin{figure}[tbh]
   \centering   
   \subfigure[Execution time: sequential.]{
     \includegraphics[width=.47\textwidth]{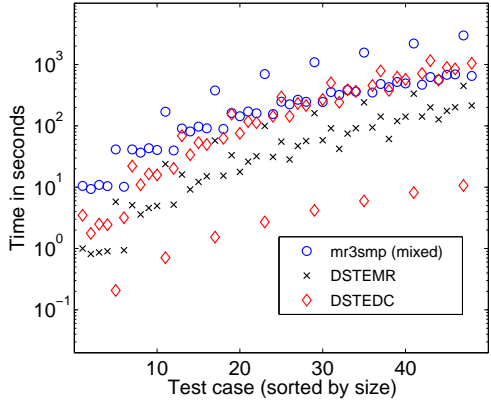}
     \label{fig:accartifial32corea}
   } \subfigure[Execution time: multi-threaded.]{
     \includegraphics[width=.47\textwidth]{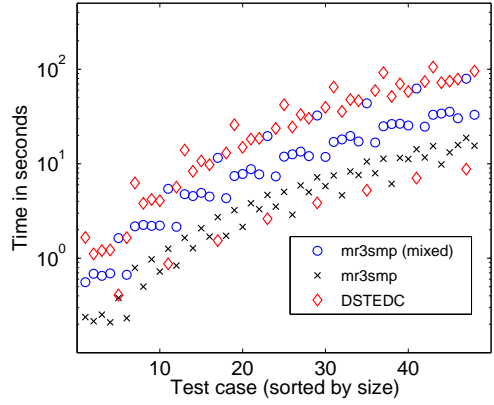}
     \label{fig:accartifial32coreb}
   }

   \subfigure[Largest residual norm.]{
     \includegraphics[width=.47\textwidth]{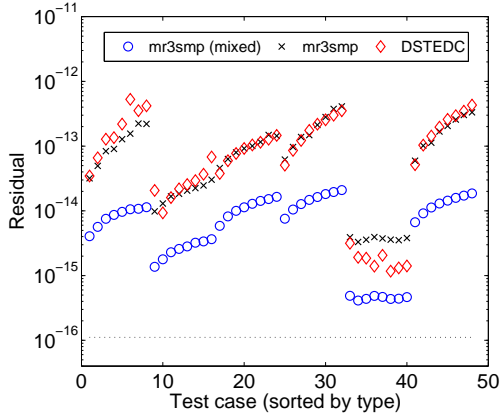}
     \label{fig:accartifial32corec}
   } \subfigure[Orthogonality.]{
     \includegraphics[width=.47\textwidth]{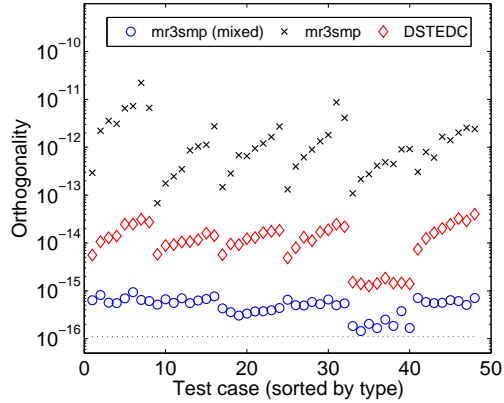}
     \label{fig:accartifial32cored}
   }
   \caption{
     Time and accuracy on {\sc Beckton}. Timings are presented in a
     logarithmic scale. The largest residual and the orthogonality are
     measured as in~\eqref{def:defresortho}. The dotted black line
     corresponds to unit round-off $\varepsilon_d$ and the accuracy of {\tt mr3smp} is
     equivalent to the one obtained by LAPACK's {\tt DSTEMR}.   
   }
   \label{fig:accartifial32core}
\end{figure}

We now turn our attention to double precision inputs/outputs, for which timings and
accuracy are presented in Fig.~\ref{fig:accartifial32core}. 
We included the results for the multi-threaded solver {\tt mr3smp}, which 
in the sequential case is just a wrapper to {\tt DSTEMR}. In general,  {\tt
  mr3smp} obtains accuracy equivalent to 
LAPACK's {\tt DSTEMR}.  

Figure~\ref{fig:accartifial32corea} shows the timings for 
sequential executions. Our mixed precision solver is slower than {\tt DSTEMR}, which is not
a surprise, as we make use of {\it software-simulated} quadruple precision
arithmetic. What might be a surprise is that, even with the use of such slow
arithmetic, for large matrices, our solver is usually as fast as {\tt
  DSTEDC}. As in the single precision case, only for matrices that allow
for substantial deflation, {\tt
  DSTEDC} is considerably faster. As Fig.~\ref{fig:accartifial32coreb}
shows, for a parallel execution, the performance difference reduces and is
expected to eventually vanish as it does already for the a regular MRRR
implementation~\cite{EleMRRR}. For matrices that do not allow for extensive
deflation, our solver is about a factor two faster than {\tt DSTEDC}. 

 While {\tt DSTEMR} accepted in 29 out of the 48 test cases representations
 that did not pass the test for relative robustness, our mixed precision
 solver found suitable representations in all cases. In fact, 
for all but the Wilkinson type matrices, we have $d_{max} = 0$ and as a
consequence: {\it no danger of failing} to find suitable representations and
{\it embarrassingly parallel} computation. 
Even for Wilkinson type matrices, $d_{max}$ was limited to one
and clustering $\rho$ was limited to $2/n$. 
For {\tt DSTEMR}, $d_{max}$ was as high as 21 and clustering $\rho$ was
about $0.7$ on average, which should 
be compared with the value of about $1.6 \cdot 10^{-4}$ for the mixed precision solver.
Therefore, we
believe that our approach is especially well-suited for highly parallel
systems. In particular, solvers for distributed-memory systems should greatly
benefit from better load-balancing and reduced communication. 

\begin{figure}[tbh]
   \centering   
   \subfigure[Execution time: multi-threaded.]{
     \includegraphics[width=.47\textwidth]{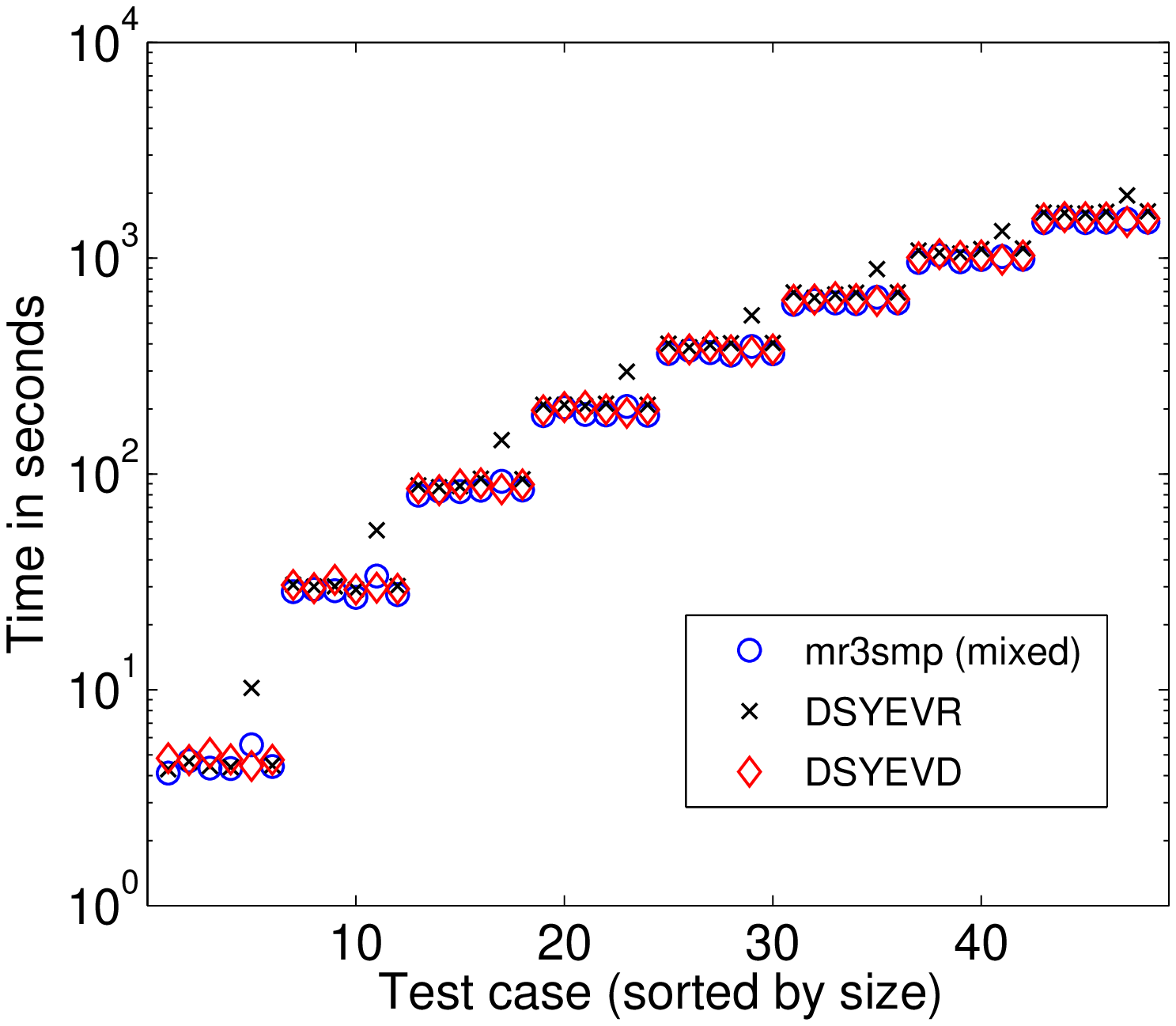}
     \label{fig:accartifial32coredensea}
   } \subfigure[Orthogonality.]{
     \includegraphics[width=.47\textwidth]{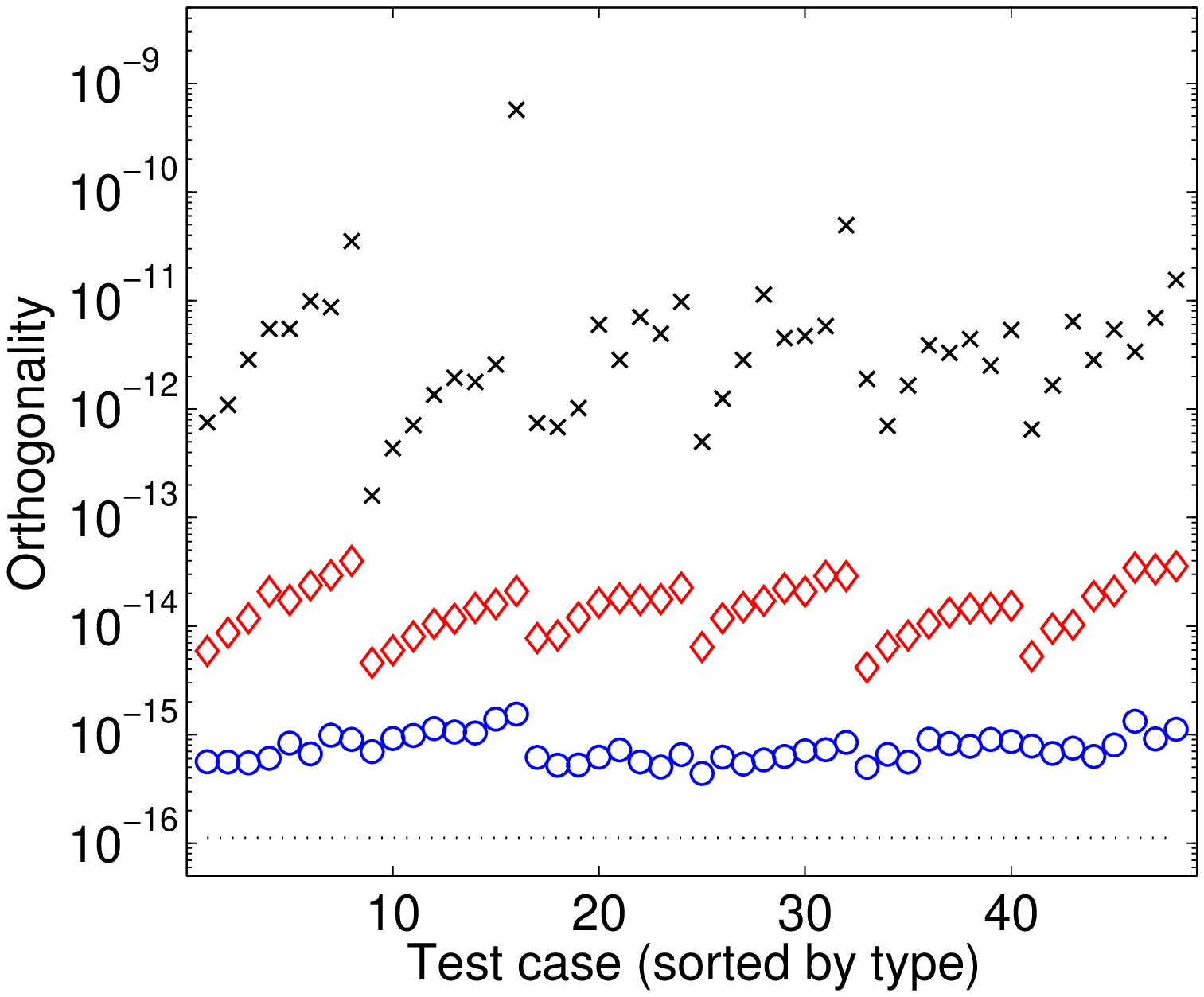}
     \label{fig:accartifial32coredenseb}
   }
   \caption{
     Time and orthogonality for computing all eigenpairs of dense, real
     symmetric matrices. Timings are presented in a
     logarithmic scale and are dominated by the reduction to tridiagonal
     form.    
   }
   \label{fig:accartifial32coredense}
\end{figure}

For single precision inputs [Figs.~\ref{fig:accartifial32coresinglea}--\subref{fig:accartifial32coresingleb}]
or in a parallel stetting [Fig.~\ref{fig:accartifial32coreb}], our tridiagonal
eigensolver is highly competitive in terms of execution time --
often faster -- compared with Divide \& Conquer and the
conventional MRRR. As a consequence, when used in context of dense Hermitian
eigenproblems, the accuracy improvement of the tridiagonal stage carry over to the
overall accuracy {\it without any penalty in terms of performance}. Such a
behavior is illustrated by
Fig.~\ref{fig:accartifial32coredense}, where we present timings and orthogonality
for dense, real symmetric input matrices. 
The inputs were generated by applying random orthogonal similarity
transformations to the tridiagonal matrices of the previous experiments:
$A = QTQ^*$, with random orthogonal matrix $Q \in \Rnn$. 
For small matrices in a sequential execution, our approach
introduces extra overhead -- see Fig.~\ref{fig:accartifial32corea}. Since 
the tridiagonal solver only requires $\order{kn}$ operations to compute $k$
eigenpairs, while the reduction to tridiagonal form requires
$\order{n^3}$ operations, for sufficiently large matrices the overhead is
completely negligible. {\em Such a statement would even more apply if the matrices were
complex-valued and/or only a subset of eigenpairs were computed, since the
reduction to tridiagonal form would carry even more
weight relative to the tridiagonal stage.} 
In a parallel execution, the mixed precision approach is competitive even
for relatively small matrices [Fig.~\ref{fig:accartifial32coredensea}]; at
the same time, the approach significantly 
improves orthogonality [Fig.~\ref{fig:accartifial32coredenseb}]. Further
experiments, including  
subset computations and complex-valued inputs, can be found in~\cite{mydiss,mixedtr}.

\section{Conclusions}

We presented a mixed precision variant of the MRRR algorithm, which addresses a
number potential weaknesses of MRRR such as ($i$) inferior
accuracy compared with the Divide \& Conquer method or the QR algorithm;
($ii$) the danger of not finding suitable representations; and ($iii$) for distributed-memory
architectures, load-balancing and
communication problems for matrices with large clustering of the
eigenvalues. Our approach provides a
new perspective: Given input/output arguments in a \binaryx\ floating
point format, we use a higher precision \binaryy\ arithmetic to
obtain the desired accuracy. As our analysis shows, the use of higher
precision provides us with freedom in setting important parameters of the
algorithm. In particular, we select these parameters to reduce the operation count,
increase robustness, and improve parallelism; at the same time, we meet
more stringent accuracy goals. Due to these changes, our mixed precision
approach is not only as accurate as % the best available methods 
the Divide \& Conquer method or the QR algorithm but -- under many
circumstances --  is also faster than  
these methods or even faster than a conventional implementation of MRRR. 

This work was mainly motivated by the results of MRRR-based eigensolvers for
dense Hermitian problems~\cite{EleMRRR}. In the context of dense
eigenproblems, the tridiagonal stage is often completely negligible
in terms of execution time: to compute $k$ eigenpairs of a tridiagonal matrix, it
only requires $\order{kn}$ operations; the reduction to
tridiagonal form requires $\order{n^3}$ operations and is the performance
bottleneck. In terms of accuracy, the 
tridiagonal stage is responsible for most of the loss of orthogonality. The
natural question was whether it is possible to improve the accuracy to the
level of the best methods without sacrificing too much performance. As
our results show, this is indeed possible. In fact, our mixed precision
solvers are even more accurate than the ones based on Divide \& Conquer and
QR, and remain as fast, or faster, than the classical MRRR. 
Finally, an important feature of the mixed precision approach is a
considerably increased robustness 
and parallel scalability. 

% \section*{Acknowledgments}
% The authors would like to thank ...

\footnotesize
\bibliographystyle{abbrv}
\bibliography{mixed}

\end{document}